\documentclass[article]{IEEEtran}

\usepackage[english]{babel}
\usepackage[latin1]{inputenc}\usepackage{enumerate}
\usepackage{amsmath}
\usepackage{amsthm}
\usepackage{amstext}
\usepackage{amssymb}
\usepackage{mathrsfs}
\usepackage{mathtools}
\usepackage{color}
\usepackage{tikz}
\usetikzlibrary{arrows.meta,calc}
\usepackage{graphicx}

\def\w{{\omega}}

\fussy

\newtheorem{theorem}{Theorem}[section]
\newtheorem{assumption}{Assumption}[section]

\newtheorem{proposition}{Proposition}[section]
\newtheorem{lemma}{Lemma}[section]
\theoremstyle{remark}
\newtheorem{remark}{Remark}[section]

\IEEEoverridecommandlockouts

\allowdisplaybreaks

\begin{document}

\title{Stochastic Control Approach to Reputation Games}

\author{Nuh Ayg\"un Dalk\i ran and Serdar Y\"uksel \thanks{
N.A. Dalk\i ran is with the Department of Economics, Bilkent University, \c{C}ankaya, Ankara, 06800, Turkey; email: dalkiran@bilkent.edu.tr. S. Y\"uksel is with the Department of Mathematics and Statistics, Queen's University, Kingston, Ontario, K7L 3N6, Canada; email: yuksel@mast.queensu.ca. This paper was presented in part at the 2nd Occasional Workshop in Economic Theory at University of Graz, the 69th European Meeting of the Econometric Society, Geneva, Switzerland, and the 11th World Congress of the Econometric Society, Montreal, Canada. This research was partially supported by the Scientific and Technological Research Council of Turkey (TUBITAK) and the Natural Sciences and Engineering Research Council of Canada (NSERC).}}

 \maketitle

\begin{abstract}
Through a stochastic control theoretic approach, we analyze reputation games where a strategic long-lived player acts in a sequential repeated game against a collection of short-lived players. The key assumption in our model is that the information of the short-lived players is nested in that of the long-lived player. This nested information structure is obtained through an appropriate monitoring structure. Under this monitoring structure, we show that, given mild assumptions, the set of Perfect Bayesian Equilibrium payoffs coincide with Markov Perfect Equilibrium payoffs, and hence a dynamic programming formulation can be obtained for the computation of equilibrium strategies of the strategic long-lived player in the discounted setup. We also consider the undiscounted average-payoff setup where we obtain an optimal equilibrium strategy of the strategic long-lived player under further technical conditions. We then use this optimal strategy in the undiscounted setup as a tool to obtain a tight upper payoff bound for the arbitrarily patient long-lived player in the discounted setup. Finally, by using measure concentration techniques, we obtain a refined lower payoff bound on the value of reputation in the discounted setup. We also study the continuity of equilibrium payoffs in the prior beliefs.
\end{abstract}

\section{Introduction}

Reputation plays an important role in long-run relationships. When one considers buying a product from a particular firm, his action (buy/not buy) depends on his belief about this firm, i.e., the firm's reputation, which he has formed based on previous experiences (of himself and of others). Many interactions among rational agents are repeated and are in the form of long-run relationships. This is why game theorists have been extensively studying the role of reputation in long-run relationships and repeated games \cite{MailathSamuelson}. By defining reputation as a conceptual as well as a mathematical quantitative variable, game theorists have been able to explain how reputation can rationalize intuitive equilibria, as in the expectation of cooperation in early rounds of a finitely repeated prisoners' dilemma \cite{krepsmilgromrobertswilson82}, and entry deterrence in the early rounds of the chain store game \cite{krepswilson82}, \cite{milgromroberts82}. 

Recently, there has been an emergence of use of tools from information and control theory in the reputations literature (see e.g., \cite{gossner11}, \cite{ekmekcigossnerwilson12}, \cite{faingold14}). Such tools have been proved to be useful in studying various bounds on the value of reputation.

In this paper, by adopting and generalizing recent results from stochastic control theory, we provide a new approach and establish refined results on reputation games. Before stating our contributions and the problem setup more explicitly, we provide a brief overview of the related literature in the following subsection.

\subsection{Related Literature}

Kreps, Milgrom, Roberts, and Wilson \cite{krepsmilgromrobertswilson82,krepswilson82,milgromroberts82} introduced the adverse selection approach to study reputations in finitely repeated games. Fudenberg and Levine \cite{fudenberglevine89, fudenberglevine92} extended this approach to infinitely repeated games and showed that a patient long-lived player facing infinitely many short-lived players can guarantee himself a payoff close to his Stackelberg payoff when there is a slight probability that the long-lived player is a commitment type who always plays the stage game Stackelberg action. When compared to the {\it folk theorem} \cite{fudenbergmaskin86, fudenberglevinemaskin94}, their results imply an intuitive expectation: the equilibria with relatively high payoffs are more likely to arise due to reputation effects. Even though the results of Fudenberg and Levine \cite{fudenberglevine89,fudenberglevine92} hold for both perfect and imperfect public monitoring, Cripps, Mailath, and Samuelson \cite{crippsmailathsamuelson04} showed that reputation effects are not sustainable in the long-run when there is imperfect public monitoring. In other words, under imperfect public monitoring it is impossible to maintain a permanent reputation for playing a strategy that does not play an equilibrium of the complete information game. There has been further literature which studies the possibility/impossibility of maintaining permanent reputations, we refer the reader to \cite{ekmekci11, ekmekcigossnerwilson12, atakanekmekci12, atakanekmekci13, atakanekmekci14, ozdogan14, liu11, faingoldsannikov11,hornerlovo09, faingold14}.

Sorin \cite{sorin99} unified and improved some of the results in reputations literature by using tools from Bayesian learning and merging due to Kalai and Lehrer \cite{KalaiLehrer1993,KalaiLehrer1994}. Gossner \cite{gossner11} utilized relative entropy (that is, information divergence or Kullback-Leibler divergence) to obtain bounds on the value of reputations; these bounds coincide in the limit (as the strategic long-lived player becomes arbitrarily patient) with the bounds provided by Fudenberg and Levine \cite{fudenberglevine89,fudenberglevine92}. 

Recently, there have been a number of related results in the information theory and control literature on real-time signaling which provide powerful structural, topological, and operational results that are in principle similar to the reputations models analyzed in the game theory literature, despite the simplifications that come about due to the fact that in these fields, the players typically have a common utility function. Furthermore, such studies typically assume finitely repeated setups, whereas we also consider here infinitely repeated setups, which require non-trivial generalizations (see e.g., \cite{Witsenhausen,WalrandVaraiya,Teneketzis,MahajanTeneketzisJSAC,YukIT2010arXiv, wood2016optimal, YukLinZeroDelay, BorkarMitterTatikonda} for various contexts but note that all of these studies except \cite{BorkarMitterTatikonda, wood2016optimal,YukLinZeroDelay} have focused on finite horizon problems). 

Using such tools from stochastic control theory and zero-delay source coding, we provide new techniques to study reputations. These techniques not only result in a number of conclusions re-affirming certain results documented in the reputations literature, but also provide new results and interpretations as we briefly discuss in the following.

{\bf Contributions of the paper.} Our findings contribute to the reputations literature by obtaining structural and computational results on the equilibrium behavior in finite-horizon, infinite-horizon, and undiscounted settings in sequential reputation games, as well as refined upper and lower bounds on the value of reputations: We analyze reputation games where a strategic long-lived player acts in a repeated sequential-move game against a collection of short-lived players each of whom plays the stage game only once but observes signals correlated with interactions of the previous short-lived players. The key assumption in our model is that the information of the short-lived players is nested in that of the long-lived player in a causal fashion. This nested information structure is obtained through an appropriate monitoring structure. Under this monitoring structure, we obtain stronger results than what currently exists in the literature in a number of directions: (i) Given mild assumptions, we show that the set of Perfect Bayesian Equilibrium payoffs coincide with the set of Markov Perfect Equilibrium payoffs. (ii) A dynamic programming formulation is obtained for the computation of equilibrium strategies of the strategic long-lived player in the discounted setup. (iii) In the undiscounted setup, under further technical conditions, we obtain an optimal strategy for the strategic long-lived player. In particular, we provide new techniques to investigate the optimality of mimicking a Stackelberg commitment type in the undiscounted setup. (iv) The optimal strategy we obtain in the undiscounted setup also lets us obtain, through an Abelian inequality, an upper payoff bound for the arbitrarily patient long-lived player\textemdash in the discounted setup. We show that this achievable upper bound is identified with a stage game Stackelberg equilibrium payoff. (v) By using measure concentration techniques, we obtain a refined lower payoff bound on the value of reputation for a fixed discount factor. This lower bound coincides with the lower bounds identified by  Fudenberg and Levine \cite{fudenberglevine92} and Gossner \cite{gossner11} as the long-lived player becomes arbitrarily patient, i.e., as the discount factor tends to 1. (vi) Finally, we establish conditions for the continuity of equilibrium payoffs in the priors.

In the next section, we present preliminaries of our model as well as two motivating examples. Section III provides our structural results leading to the equivalence of Perfect Bayesian Equilibrium payoffs and Markov Perfect Equilibrium payoffs in the discounted setup. Section IV provides results characterizing the optimal behavior of the long-lived player for the undiscounted setup, which lead us to an upper bound for the equilibrium payoffs in the discounted setup when the long-lived player becomes arbitrarily patient.  Section V studies the continuity problem in the priors. Section VI provides, through an explicit measure concentration analysis, a refined lower bound for the equilibrium payoffs of the strategic long-lived player in the discounted setup. 

\section{The Model}

A long-lived player (Player 1) plays a repeated stage game with a sequence of different short-lived players (each of whom is referred to as Player 2).  

{\bf The stage game.} The stage game is a \textit{sequential-move} game: Player 1 moves first; when action $a^1$ is chosen by Player 1 in the stage game; a \textbf{public} signal $s^2 \in \mathbb{S}^2$ is observed by Player 2 which is drawn according to the probability distribution $\rho^2(.|a^1) \in \Delta(\mathbb{S}^2)$. Player 2, observing this public signal (and all preceding public signals), moves second. At the end of the stage game, Player 1 observes a \textbf{private} signal $s^1 \in \mathbb{S}^1$ which depends on actions of both players in the stage game and is drawn according to the probability distribution $\rho^1(.|(a^1, a^2))$.  That is, the stage game can be considered as a Stackelberg  game with imperfect monitoring where Player 1 is the leader and Player 2 is the follower. Action sets of Player 1 and Player 2 in the stage game are assumed to be finite and denoted by $\mathbb{A}^1$ and $\mathbb{A}^2$, respectively. We also assume that the set of Player 1's all possible private signals, denoted by $\mathbb{S}^1$,  and the set of (Player 2s') all possible public signals, denoted by $\mathbb{S}^2$,  are finite.

{\bf The information structure.} There is incomplete information regarding the type of the long-lived Player 1. Player 1 can either be a strategic type (or normal type), denoted by $\w^n$, or one of finitely many simple commitment types. Each of these commitment types is committed to simply playing the same action $\hat{\omega} \in \Delta(\mathbb{A}^1)$ at every stage of the repeated game\textemdash independent of the history of the play.\footnote{$\Delta(\mathbb{A}^i)$ denotes the set of all probability measures on $\mathbb{A}^i$ for both $i=1,2$. That is, the commitment types can be committed to playing mixed  stage-game actions as well. We would like to also note here that simple commitment types assumption is a standard assumption in reputation games.} The set of all possible commitment types of Player 1 is given by $\hat{\Omega}$.  Therefore, the set of all possible types of Player 1 can be denoted as  $\Omega=\{\w^n\} \cup \hat{\Omega}$. The type of Player 1 is determined once and for all at the beginning of the game according to a {\it common knowledge} and {\it full-support} prior $\mu_0 \in \Delta(\Omega)$. Only Player 1 is informed of his type, i.e., Player 1's type is his private information.

We note that there is a \textbf{nested information structure} in the repeated game in the following sense: The signals observed by Player 2s are public, and hence available to all subsequent players, whereas Player 1's signals are his private information. Therefore, the information of Player 2 at time $t-1$ is a subset of the information of Player 1 at time $t$. 
Formally, a generic history for Player 2 at time $t-1$ and a generic history for Player 1 at time $t$ are given as follows:
\begin{eqnarray}
h^2_{t-1}&=&(s^2_0,s^2_1,\cdots,s^2_{t-1}) \in  H^2_{t-1} \\
h^1_{t}&=& (a^1_0, s^1_0,s^2_0,\cdots,a^1_{t-1},s^1_{t-1},s^2_{t-1}) \in H^1_t \label{P1hist}
\end{eqnarray}
where   $ H^2_{t-1}\label{P2hist} := {(\mathbb{S}^2)}^t $ and $H^1_t \label{P1hist} := {(\mathbb{A}^1 \times \mathbb{S}^1 \times \mathbb{S}^2)}^t$.

 That is, each Player 2 observes, before he acts, a finite sequence of public signals which are correlated with Player 1's action in each of his interaction with preceding Player 2s. On the other hand, Player 1 observes not only these public signals, but also a sequence of private signals for each particular interaction that happened in the past, and his actions in the previous periods\textemdash but not necessarily the actions of preceding Player 2s.\footnote{Note that Player 1 gets to observe the realizations of his earlier possibly mixed actions.}

We note also that having such a monitoring structure is not a strong assumption. In particular,  it is weaker than the information structure in Fudenberg and Levine \cite{fudenberglevine92} where it is assumed that only the same sequence of public signals are observable by the long-lived and short-lived players, i.e., there is only public monitoring. Yet, it is stronger than the information structure in Gossner \cite{gossner11} which allows private monitoring for both the long-lived and the short-lived players. 

The stage game payoff function of the strategic (or normal) type long-lived Player 1 is given by $u^1$, and each short-lived Player 2's payoff function is given by $u^2$, where $u^i :\mathbb{A}^1 \times \mathbb{A}^2 \rightarrow \mathbb{R}$. The set of all possible histories for Player 2 of stage $t$ is $H^2_{t} = H^2_{t-1} \times \mathbb{S}^2$  where $H^2_{t-1} = {(\mathbb{S}^2)}^t$. On the other hand, the set of all possible histories observable by the long-lived Player 1 prior to stage $t$ is $H^1_{t} ={(\mathbb{A}^1 \times \mathbb{S}^1 \times \mathbb{S}^2)}^t $. It is assumed that $H^1_{0} := \emptyset$ and $H^2_{0} := \emptyset$, which is the usual convention. Let $\mathcal{H}^1 = \bigcup_{t\geq 0} H^1_{t} $ be the set of all possible histories of the long-lived Player 1.

A (behavioral) strategy for Player 1 is a map:
\begin{equation*}
\sigma^1 :\Omega \times \mathcal{H}^1\rightarrow \Delta(\mathbb{A}^1).
\end{equation*}
that satisfies $\sigma^1(\hat{\w},h^1_{t-1}) = \hat{\w}$ for any $\hat{\w} \in \hat{\Omega}$ and for every $h^1_{t-1} \in H^1_{t-1}$, since commitment types are required to play the corresponding (fixed) action of the stage game independent of the history. The set of all strategies for Player 1 is denoted by $\Sigma^1$, i.e., $\Sigma^1$ is the set of all functions from ${\Omega \times \mathcal{H}^1}$ to $\Delta(\mathbb{A}^1)$.

A strategy for Player 2 of stage $t$ is a map:
\begin{equation*}
\sigma^2_{t}: H^2_{t-1} \times \mathbb{S}^2 \rightarrow \Delta(\mathbb{A}^2).	
\end{equation*}

We let $\Sigma^2_t$ be the set of all such strategies and let $\Sigma^2= \Pi_{t \geq 0}\Sigma^2_t$ denote the set of
all sequences of all such strategies. A history (or path) $h_t$ of length $t$ is an element of $\Omega \times (\mathbb{A}^1\times \mathbb{A}^2 \times \mathbb{S}^1 \times \mathbb{S}^2)^t$ describing Player 1's type, actions, and signals realized up to stage $t$. By standard arguments (e.g., Ionescu-Tulcea Theorem \cite{HernandezLermaMCP}), a strategy profile $\sigma =(\sigma^1, \sigma^2) \in \Sigma^1 \times \Sigma^2$ induces a unique probability distribution $P_{\sigma}$ over the set of all paths of play $H^\infty = \Omega \times (\mathbb{A}^1 \times \mathbb{A}^2 \times \mathbb{S}^1 \times \mathbb{S}^2)^{\mathbb{Z}_+}$ endowed with the product $\sigma$-algebra. We let $a_t =(a^1_{t}, a^2_t)$ represent the action profile realized at stage $t$ and let $s_t =
(s^1_{t}, s^2_{t})$ denote the signal profile realized at stage $t$. Given $\omega \in \Omega$, $P_{\w, \sigma}(.)= P_{\sigma}(.|\w)$ represents the probability distribution over all paths of play conditional on Player 1 being type $\w$. Player 1's discount factor is assumed to be $\delta \in (0,1)$ and hence, the expected discounted average payoff to the strategic (normal type) long-lived
Player 1 is given by
\begin{equation*}
\pi_1(\sigma)=\mathbb{E}_{P_{\w^n,\sigma}} (1-\delta) \sum_{t\geq 0} \delta^{t} u^1 (a_t).
\end{equation*}

In all of our results except Lemma \ref{witRealCoding}, we will assume that Player 2s are Bayesian rational.\footnote{A Bayesian rational Player 2 tries to maximize his expected payoff after updating his beliefs according to the Bayes' rule whenever possible. We also note that Lemma \ref{witRealCoding} does not require Bayesian rationality and holds for non-Bayesian Player 2s who might underreact or overreact to new (or recent) information as in \cite{epsteinnoorsandroni10} as well.} Hence, we will restrict attention to Perfect Bayesian Equilibrium: In any such equilibrium, the strategic Player 1 maximizes his expected discounted average payoff given that the short-lived players play a best response to their expectations according to their updated beliefs (This will  be appropriately modified when we consider the undiscounted setup). Each Player 2, playing the stage game only once, will be best-responding to his expectation according to his beliefs which are updated according to the Bayes' Rule.

A strategy of Player 2s, $\sigma^{2}$, is a best response to $\sigma^{1}$ if, for all $t$,
\begin{eqnarray*}
\mathbb{E}_{P_{\sigma }}[u^{2}(a^1_{t},a^2_{t})|s^2_{[0,t]}]\geq \mathbb{E}_{P_{\sigma }}[u^{2}(a^1_{t},a^2)|s^2_{[0,t]}] \\ \text{for all } a^2\in A^2 \text{    }\ (P_{\sigma}-a.s.)
\end{eqnarray*}
where $s^2_{[0,t]} =(s^2_0,s^2_1,\cdots,s^2_{t})$ denotes the information available to Player 2 at time $t$.

%Next, we provide two examples to motivate our model.

\subsection{Motivating Example I: The Product Choice Game}

Our first example is a simple product choice game which describes how a strategic player can build up reputation: There is a (long-lived) firm (Player 1) who faces an infinite sequence of different consumers (Player 2s) with identical preferences. There are two actions available to the firm: $A_{1}=\{H,L\}$, where $H$ and $L$ denote exerting \emph{high-effort} and \emph{low-effort} in the production of its output, respectively. Each consumer also has two possible actions: buying a \emph{high-priced} product, $(h)$, or a \emph{low-priced} product, $(l)$, i.e., $A_{2}=\{h,l\}$.  Each consumer prefers a high-priced product if the firm exerted high effort and a low-priced product if the firm exerted low effort. The firm is willing to commit to high effort only if the consumers purchase the high-priced product, i.e., the firm's (pure) Stackelberg action \textemdash in the stage game\textemdash\ is exerting high level of effort. Therefore, if the level of effort of the firm were observable, each consumer would best reply to the Stackelberg action by buying a high priced product. However, the choice of effort level of the firm is not observable before consumers choose the product. Furthermore, exerting high effort is costly, and hence, for each type of product, the firm prefers to exert low effort rather than high effort. That is, there is a moral hazard problem.

The corresponding stage game and the preferences regarding the stage game can be illustrated as follows:

\begin{center}
\begin{tikzpicture}
  \newcommand{\XS}{3}
  \newcommand{\YS}{2}
  \coordinate (PODN1) at (0,0);
  \coordinate (PTDN1) at ($ (PODN1) + (\XS,\YS) $);
  \coordinate (PTDN2) at ($ (PODN1) + (\XS,-\YS) $);
  \coordinate (PTDN3) at ($ (PODN1) + (\XS,0.1*\YS - 0.9*\YS) $);
  \coordinate (PTDN4) at ($ (PODN1) + (\XS,0.30*\YS - 0.70*\YS) $);
  \coordinate (PTDN5) at ($ (PODN1) + (\XS,0.50*\YS - 0.50*\YS) $);
  \coordinate (PTDN6) at ($ (PODN1) + (\XS,0.70*\YS - 0.30*\YS) $);
  \coordinate (PTDN7) at ($ (PODN1) + (\XS,0.9*\YS - 0.1*\YS) $);
  \coordinate (TN1) at ($ (PTDN2) + (\XS,\YS) $);
  \coordinate (TN2) at ($ (PTDN2) + (\XS,-\YS) $);
  \coordinate (TN3) at ($ (PTDN2) + (\XS,0.1*\YS - 0.9*\YS)  $ );
  \coordinate (TN4) at ($ (PTDN4) + (\XS, 0.1\YS - 0.7*\YS)  $ );
  \coordinate (TN5) at ($ (PTDN2) + (\XS,0.5*\YS - 0.5*\YS)  $ );
  \coordinate (TN6) at ($ (PTDN2) + (\XS,0.7*\YS - 0.3*\YS)  $ );
  \coordinate (TN7) at ($ (PTDN2) + (\XS,0.9*\YS - 0.1*\YS)  $ );
  \coordinate (TN8) at ($ (PTDN5) + (\XS,\YS)  $ );
  \coordinate (TN9) at ($ (PTDN7) + (\XS, 0.75*\YS - 0.25*\YS)  $ );
  \coordinate (TN10) at ($ (PTDN7) + (\XS,0.25*\YS - 0.75*\YS)  $ ); 
  \coordinate (TN11) at ($ (PTDN3) + (\XS,0.75*\YS - 0.25*\YS)  $ );
  \coordinate (TN12) at ($ (PTDN3) + (\XS,0.25*\YS - 0.75*\YS)  $ );

   \draw [dashed] (PTDN3) -- node [left] {$P_2$}(PTDN7);
  \draw (PODN1)  node [left] {$P_1$} --  node [above right, sloped]{$H$} (PTDN7);
  \draw  (PODN1) -- node [below right, sloped]{$L$} (PTDN3);

\node at  ($(TN9) + (0.4,0)$) {$(2,3)$};
\node at  ($(TN10) + (0.4,0)$) {$(0,2)$};
\node at ($(TN11) + (0.4,0)$) {$(3,0)$};
\node at  ($(TN12) + (0.4,0)$) {$(1,1)$};

  \draw  (PTDN7) --  node [above right, sloped]{$h$} (TN9);
  \draw (PTDN7) --  node [below right, sloped]{$l$}(TN10);
  \draw  (PTDN3) --  node [above right, sloped]{$h$} (TN11);
  \draw (PTDN3) --  node [below right, sloped]{$l$}(TN12);

  \coordinate (UL) at (current bounding box.north west);
  \coordinate (LL) at (current bounding box.south west);
  \coordinate (UR) at (current bounding box.north east);
  \coordinate (LR) at (current bounding box.south east);

  \draw
  [{<[scale=1.8,length=2.25,width=2.25]}-{>[scale=1.8,length=2.25,width=2.25]}]
%  [{<[scale=2.5,length=3,width=3]}-{>[scale=2.5,length=3,width=3]}]
  ([yshift = -0.5cm]LL) -- node [below] {  \textbf{Figure 1}: The illustration of the stage game }([yshift =
  -0.5cm]LR) ;
\end{tikzpicture}
\end{center}

\begin{center}
\quad \quad \quad \quad
\begin{tabular}{ccc}
& $h$ & $l$ \\ \cline{2-3}
$H$ & \multicolumn{1}{|c}{$2,3$} & \multicolumn{1}{|c|}{$0,2$} \\
\cline{2-3}
$L$ & \multicolumn{1}{|c}{$3,0$} & \multicolumn{1}{|c|}{$1,1$} \\
\cline{2-3}
\end{tabular}
\newline
\end{center}

Note that since the stage game is a sequential-move game where actions are not observable, it is strategically equivalent to a simultaneous-move game represented by the corresponding payoff matrix, which is given above. Furthermore, there is a \emph{unique} Nash equilibrium of this stage-game, and in this equilibrium the firm (the row player) plays $L$ (exerts \emph{low effort}) and the consumer (the column player) plays $l$ (buying a \emph{low-priced} product). 

Suppose that there is a small but positive probability $p_{0}>0$ that the firm is an honorable firm who always exerts \emph{high effort}. That is, with $p_{0}>0$ probability, Player 1 is a \emph{commitment type} who plays $H$ at every period of the repeated game\textemdash independent of the history.  Suppose further that each consumer can observe all the outcomes of the previous play. Yet, before he acts, the consumer cannot observe the effort level of the firm in his own period of play.

Consider now a strategic (non-commitment or normal type) firm who has a discount factor $\delta <1$: Can the firm build up a reputation that he is (or acts as if he is) an honorable firm? The answer to this question is ``Yes''\textemdash when he is patient enough. 

To see this, observe that a rational consumer (Player 2) would play $h$ only if he anticipates that the firm (Player 1) plays $H$ with a probability of at least $\frac{1}{2}$. Let $p_t$ be the posterior belief that Player 1 is a commitment type after observing some public history $h_t$. Suppose Player 2 of period $t+1$ observes $(H,l)$ as the outcome of the preceding period $t$. This means the probability that Player 2 of period $t$ anticipated for $H$ was less than (or equal) to $\frac{1}{2}$. This probability is $p_t+(1-p_t)\sigma^1(\omega^n, h_t)(H)$  where $\sigma^1(\omega^n, h_t)(H)$ is the probability that the strategic (or normal) type Player 1 assigns to playing $H$ at period $t$ after observing $h_t$. Therefore, we have $p_t+(1-p_t)\sigma^1(\omega^n, h_t)(H)\leq \frac{1}{2}$. But, this implies that the posterior belief of Player 2 of period $t+1$ that Player 1 is a commitment type \textemdash after observing $(H,l)$\textemdash\  will be $p_{t+1}= \frac{p_t}{p_t+(1-p_t)\sigma^1(\omega^n, h_t)(H)}\geq 2 p_t$. This means every time the strategic player plays $H$, he doubles his reputation, i.e., the belief that he is a commitment type doubles. Therefore, mimicking the commitment type finitely many rounds, the firm can increase the belief that he is an honorable firm (a commitment type) with more than probability $\frac{1}{2}$. In such a case, the short lived consumers (Player 2s) will start best replying by buying high-priced products. If the firm is patient enough \textemdash when $\delta$ is high\textemdash\ payoffs from those finitely many periods will be negligible. Furthermore, as $\delta \rightarrow 1 $, one can show that the strategic Player 1 can guarantee himself a discounted average payoff arbitrarily close to $2$\textemdash which is his payoff under his  (pure) Stackelberg action. 

\subsection{Motivating Example II: A Consultant with Reputational Concerns under Moral Hazard}

Our second example presents finer details regarding the nested information structure: 
A consultant is to advise different firms in different projects. In each of these projects, a supervisor from the particular firm is to inspect the consultant regarding his effort during the particular project. The consultant can either exert a (H)igh level of effort or a (L)ow level of effort while working on the project.

The effort of the consultant is not directly observable to the supervisor. Yet, after the consultant chooses his effort level, the supervisor gets to observe a public signal $s^2 \in \{h,l\}$ which is correlated with the effort level of the consultant according to the probability distribution $\rho^2(h|H)=\rho^2(l|L)= p >\frac{1}{2}$.

Observing this public signal, the supervisor recommends to the upper administration to give the consultant a (B)onus or (N)ot. 

The supervisor prefers to recommend a (B)onus when the consultant works hard (exerts (H)igh effort) and  (N)ot to recommend a bonus when the consultant shirks (exerts (L)ow effort). For the consultant exerting a high level of effort is costly. Therefore, the stage game and the preferences regarding the stage game can be illustrated as follows:\footnote{Note that the stage game is a sequential-move game, the payoffs are summarized in a payoff matrix just for illustrative purposes.}

\begin{center}
\begin{tikzpicture}
  \newcommand{\XS}{3}
  \newcommand{\YS}{2}
  \coordinate (PODN1) at (0,0);
  \coordinate (PTDN1) at ($ (PODN1) + (\XS,\YS) $);
  \coordinate (PTDN2) at ($ (PODN1) + (\XS,-\YS) $);
  \coordinate (PTDN3) at ($ (PODN1) + (\XS,0.10*\YS - 0.90*\YS) $);
  \coordinate (PTDN4) at ($ (PODN1) + (\XS,0.30*\YS - 0.70*\YS) $);
  \coordinate (PTDN5) at ($ (PODN1) + (\XS,0.50*\YS - 0.50*\YS) $);
  \coordinate (PTDN6) at ($ (PODN1) + (\XS,0.70*\YS - 0.30*\YS) $);
  \coordinate (PTDN7) at ($ (PODN1) + (\XS,0.90*\YS - 0.10*\YS) $);
  \coordinate (TN1) at ($ (PTDN2) + (\XS,\YS) $);
  \coordinate (TN2) at ($ (PTDN2) + (\XS,-\YS) $);
  \coordinate (TN3) at ($ (PTDN2) + (\XS,0.1*\YS - 0.9*\YS)  $ );
  \coordinate (TN4) at ($ (PTDN4) + (\XS, 0.1\YS - 0.7*\YS)  $ );
  \coordinate (TN5) at ($ (PTDN2) + (\XS,0.5*\YS - 0.5*\YS)  $ );
  \coordinate (TN6) at ($ (PTDN2) + (\XS,0.7*\YS - 0.3*\YS)  $ );
  \coordinate (TN7) at ($ (PTDN2) + (\XS,0.9*\YS - 0.1*\YS)  $ );
  \coordinate (TN8) at ($ (PTDN5) + (\XS,\YS)  $ );

   \draw [dashed] (PTDN1) -- (PTDN2);
  \draw (PODN1)  node [left] {$P_1$} --  node [above right, sloped]{$H$}(PTDN1);
  \draw  (PODN1) -- node [below right, sloped]{$L$} (PTDN2);
  \draw  (PTDN5) node[left]{$s^2 \in\{h,l\}$} node [right] {$P_2$} --  node [below right, sloped]{$N$} (TN4);
  \draw (PTDN5) --  node [above right, sloped]{$B$}(TN8);

  \coordinate (UL) at (current bounding box.north west);
  \coordinate (LL) at (current bounding box.south west);
  \coordinate (UR) at (current bounding box.north east);
  \coordinate (LR) at (current bounding box.south east);

  \draw
  [{<[scale=2.0,length=2.5,width=2.5]}-{>[scale=2.0,length=2.5,width=2.5]}]
  ([yshift = -0.5cm]LL) -- node [below] {  \textbf{Figure 2}: The illustration of the stage game }([yshift =
  -0.5cm]LR) ;

\end{tikzpicture}

\end{center}

\begin{center}
\quad \quad \quad
\begin{tabular}{ccc}
& $B$ & $N$ \\ \cline{2-3}
$H$ & \multicolumn{1}{|c}{$1,1$} & \multicolumn{1}{|c|}{$-1,-1$} \\
\cline{2-3}
$L$ & \multicolumn{1}{|c}{$2,-2$} & \multicolumn{1}{|c|}{$0,0$} \\
\cline{2-3}
\end{tabular}
\newline
\end{center}

It is commonly known that there is a positive probability $p_{0}>0$ with which the consultant is an honorable consultant  who always exerts (H)igh level of effort. That is, with $p_{0}>0$ probability the consultant is a \emph{commitment type}  who plays $H$ at every period of the repeated game independent of the history.

Consider the incentives of a strategic (non-commitment or normal type) consultant:  Does such a consultant have an incentive to build a reputation by exerting high level of effort, if the game is repeated only finitely many times? What kind of equilibrium behavior would one expect from such a consultant  if the game  is repeated infinitely many times with discounting for {\it a fixed discount factor}? For example, if he is building a reputation, how often does he shirk (exert (L)ow level of effort)?  Does there exist reputation cycles, i.e., does the consultant build a reputation by exerting high effort for a while and then milks it by exerting low effort until his reputation level falls under a particular threshold? What happens when the consultant becomes arbitrarily patient, i.e., his discount factor tends to 1? What can we say about the consultant's optimal reputation building strategy when he does not discount the future but rather cares about his undiscounted average payoff?

The aim of this paper is to provide tractable techniques to answer similar questions in settings where agents have reputational concerns in repeated game setups described in our model.

\section{Optimal Strategies and Equilibrium Behavior}

Our first set of results will be regarding the optimal strategies of the strategic long-lived Player 1.

Briefly, since each Player 2 plays the stage game only once, we show that when the information of Player 2 is nested in that of Player 1, under a plausible assumption to be noted the strategic long-lived Player 1 can, without any loss in payoff performance, formulate his strategy as a controlled Markovian system optimization, and thus through dynamic programming. The discounted nature of the optimization problem then leads to the existence of a stationary solution. This implies that for any Perfect Bayesian Equilibrium, there exists a payoff-equivalent stationary Markov Perfect Equilibrium. Hence, we conclude that the Perfect Bayesian Equilibrium payoff set and Markov Perfect Equilibrium payoff set of the strategic long-lived Player 1 coincide with each other.

Below, we provide three results on optimal strategies following steps parallel to  \cite{YukselBasarBook} which builds on Witsenhausen \cite{Witsenhausen}, Walrand and Varaiya \cite{WalrandVaraiya}, Teneketzis \cite{Teneketzis}, and \cite{YukIT2010arXiv}. These structural results on optimal strategies will be the key for the following Markov chain construction as well as Theorem \ref{pbempef} and Theorem \ref{pbempe}. 

\subsection{Optimal Strategies: Finite Horizon}\label{FiniteHorizon}

We first consider the finitely repeated game setup where the stage game is to be repeated $T \in \mathbb{N}$ times. In such a case, the strategic long-lived Player 1 is to maximize $\pi_1(\sigma)$ given by
\begin{equation*}
\pi_1(\sigma)=\mathbb{E}_{P_{\w^n,\sigma}} (1-\delta) \sum_{t=0}^{T-1} \delta^{t} u^1 (a_t).
\end{equation*}

Our first result, Lemma \ref{witRealCoding}, shows that, given any fixed sequence of strategies of the short-lived Player 2s,  any optimal strategy of the strategic long-lived Player 1 can be replaced, without any loss in payoff performance, by another optimal strategy which only depends on the (public) information of Player 2s. More specifically, we show that for any private strategy of the long-lived Player 1 against an arbitrary sequence of strategies of Player 2s, there exists a public strategy  of the long-lived Player 1 against the very same sequence of strategies of Player 2s  which gives the strategic long-lived player a better payoff.\footnote{A public strategy is a strategy that uses only public information that is available to all the players. On the other hand,  a strategy that is based on private information of a player is referred to as a private strategy. In particular, any strategy of Player 1 that is based on  $s^1_t$ for some $t$ is a private strategy. }

To the best of our knowledge, this is a new result in the repeated games literature. What is different here from similar results in the repeated games literature is that this is true even when Player 2s strategies are non-Bayesian.\footnote{A relevant result appears in \cite{fudenberglevine94}, which shows that sequential equilibrium payoffs and perfect public equilibrium payoffs coincide (See the Appendix B of \cite{fudenberglevine94}) in a similar infinitely repeated game setup.}

Before we state Lemma \ref{witRealCoding}, we note here that the signal $s^2_t$ that will be available to short-lived Player 2s after round $t$ only depends on the action of the long-lived Player 1 at round $t$ and that the following holds for all $t \geq 1$.
\begin{equation}\label{ChannelMemoryless}
P_\sigma(s^2_t| a^1_t; a^1_{t'}, a^2_{t'}, t' \leq t-1) = P_\sigma(s^2_t|a^1_t).
\end{equation}

Observation (\ref{ChannelMemoryless}) plays an important role in the proof of our first result:\\

\begin{lemma}\label{witRealCoding}
In the finitely repeated setup,  given any sequence of strategies of short-lived Player 2s, for any (private) strategy of the strategic long-lived  Player 1, there exists a (public) strategy that only conditions on  $\{s^2_0,s^2_1,\cdots, s^2_{t-1}\}$ which yields the strategic long-lived Player 1 a better payoff against the given sequence of strategies of Player 2s.
\end{lemma}

\noindent{\bf Proof.} See the Appendix. \\

A brief word of caution is in order. The structural results of the type Lemma \ref{witRealCoding}, while extremely useful in team theory and zero-delay source coding \cite{YukselBasarBook}, do not always apply to generic games unless one further restricts the setup. In particular, a generic (Nash) equilibrium may be lost once one alters the strategy structure of one of the players, while keeping the other one fixed (in team problems, the parties can agree to have a better performing team policy even if it is not a strict equilibrium). However, we consider the Perfect Bayesian Equilibrium concept here which is of a leader-follower type (i.e., {\it Stackelberg in the policy space}): {Perfect Bayesian Equilibrium requires sequential rationality--and hence eliminates non-credible threats.} That is, Player 2s respond in a Bayesian fashion to Player 1 who in turn is aware of Player 2s commitment to this policy. This subtle difference is crucial also in signaling games; the features that distinguish Nash equilibria (as in the classical setup studied in Crawford and Sobel \cite{SignalingGames}) from Stackelberg equilibria in signaling games are discussed in detail in \cite[Section 2]{saritacs2017quadratic}. 

%\sy{What is the solution concept in Gossner/Sorin etc.? Let's include a very brief discussion on this, since theirs is different evidently.}

Lemma \ref{witRealCoding} implies that any private information of Player 1 is statistically irrelevant for optimal strategies: for any private strategy of the long-lived Player 1, there exists a public strategy which performs at least as good as the original one against a given sequence of strategies of Player 2s. That is, in the finitely repeated setup, the long-lived Player 1 can depend his strategy only on the public information and his type without any loss in payoff performance. We would like to note here once again that Lemma \ref{witRealCoding} above holds for any sequence of strategies of Player 2s, even non-Bayesian ones.

On the other hand, when Player 2s are Bayesian rational, as is the norm in repeated games, we obtain a more refined structural result which we state below as Lemma \ref{walrandType0}. As mentioned before, in a Perfect Bayesian Equilibrium the short-lived Player 2 at time $t$, playing the stage game only once, seeks to maximize $ \sum_{a^1} P_{\sigma}(a^1_t=a^1|s^2_{[0,t]}) u^2(a^1,a^2)$.
However, it may be that his best response set, i.e., the maximizing action set $\mathrm{arg}\max(\sum_{a^1} P_{\sigma}(a^1_t=a^1|s^2_{[0,t]}) u^2(a^1,a^2))$, may not be unique.

To avoid such set-valued correspondence dynamics, we consider the following assumption, which requires that the best response of each Player 2 is essentially unique: Note that any strategy for Player 2 of time $t$ who chooses
\[\mathrm{arg}\max\big(\sum_{a^1} P_\sigma(a^1_t=a^1|s^2_{[0,t]}) u^2(a^1,a^2)\big)\] in a measurable fashion does not have to be continuous in the conditional probability $\kappa(\cdot) = P_\sigma(a^1_t = \cdot |s^2_{[0,t]})$, since such a strategy partitions (or quantizes) the set of probability measures on $\mathbb{A}^1$. The set of $\kappa$ which borders these partitions is a subset of the set of probability measures ${\cal B}_e = \cup_{k,m \in \mathbb{A}^2} {\cal B}^{k,m}$, where for any pair $k,m \in \mathbb{A}^2$, the belief set $ {\cal B}^{k,m}$ is defined as
\begin{equation}
\resizebox{0.446 \textwidth}{!}{${\cal B}^{k,m} =\bigg\{ \kappa \in  \Delta(\mathbb{A}^1): \sum_{a^1 \in \mathbb{A}^1} \kappa(a^1) u^2(a^1,k) = \sum_{a^1 \in \mathbb{A}^1} \kappa(a^1) u^2(a^1,m)\bigg\}.$}
\end{equation}
These are the sets of probability measures where Player 2 is indifferent between multiple actions.\footnote{In particular, in both of our motivating examples, the set $B_e$ is the singleton probability measure $\{(\frac{1}{2},\frac{1}{2})\}$. To see this, it is enough to consider the corresponding payoff matrix for each of the motivating examples. One can verify that in both of the motivating examples, Player 2 becomes indifferent only when Player 1 randomizes between H and L with $\frac{1}{2}$ probability.}

\begin{assumption}\label{EssUnique}
Either of the following holds:
\begin{itemize}
	\item[(i)] The prior measure and the probability space is so that $P_\sigma\bigg( P_\sigma(a^1_t = \cdot |s^2_{[0,t]}) \in {\cal B}_e \bigg) = 0$ for all $t \geq 0$. In particular, Player 2s have a unique best response so that the set of discontinuity, ${\cal B}_e$, is never visited (with probability $1$).
	\item[(ii)] Whenever Player 2s are indifferent between multiple actions they choose the action that is better for Player 1.
\end{itemize}
%\sy{All we need in (ii) is that Player 1 knows what Player 2 will do; should we give a more general condition? e.g., we can have the worse action too.}-- NAD: I  think there is no need. But, it is up to you. 
\end{assumption}
The following remarks are on Assumption \ref{EssUnique}.
\begin{remark} 
\begin{itemize}\
	\item[(i)] In the classical reputations literature, a standard result is that under mild conditions Bayesian rational short-lived players can be {\it surprised} at most finitely many times, e.g., \cite[Theorem 4.1]{fudenberglevine92}, \cite[Lemma 2.4]{sorin99}, implying that the jumps in the corresponding belief dynamics of Player 2s will be bounded away from zero in a transient phase until the optimal responses of Player 2s converge to a fixed action. In such cases, the payoff structure can be designed so that the set of discontinuity, ${\cal B}_e$, is visited with 0 probability, and hence Assumption \ref{EssUnique}(i) holds. %This assumption will also be utilized separately in the analysis with regard to the undiscounted setup later in the paper.
	\item[(ii)] Assumption \ref{EssUnique}(ii) is a standard assumption in the contract theory literature. In a principal-agent model, whenever an agent is indifferent between two actions he chooses the action that is better for the principal, e.g., when an incentive compatibility condition binds so that the agent is indifferent between exerting a high level of effort and exerting a low level effort, then the agent chooses to exert the high level of effort (see \cite{BoltonDewatripont} for further details). Assumption \ref{EssUnique}(ii) trivially holds also when the stage game payoff functions are identical for both players (as in {\it team} setups) or are aligned (as in a {\it potential game}).  

\end{itemize}
\end{remark} 

\begin{lemma}\label{walrandType0}
	In the finitely repeated setup, under Assumption \ref{EssUnique}, given any arbitrary sequence of strategies of Bayesian rational short-lived Player 2s, for any (private) strategy of the strategic long-lived  Player 1, there exists a (public) strategy that only conditions on  $P_\sigma(\omega | s^2_{[0,t-1]}) \in \Delta(\Omega)$ and $t$  which yields the strategic long-lived Player 1 a better payoff against the given sequence of strategies of Player 2s.
\end{lemma}
\noindent{\bf Proof.} See the Appendix.

\subsection{Controlled Markov Chain Construction}\label{MarkovChain}

The proof of Lemma \ref{walrandType0} reveals the construction of a controlled Markov chain. Building on this proof, we will explicitly construct the dynamic programming problem as a controlled Markov chain optimization problem (that is, a {\it Markov Decision Process}). Under Assumption \ref{EssUnique}, given any sequence of strategies of Bayesian rational Player 2s, the solution to this optimization problem characterizes the equilibrium behavior of the strategic long-lived player in an associated Markov Perfect Equilibrium. The state space, the action set, the transition kernel, and the per-stage reward function of the controlled Markov chain mentioned above are given as follows:

\begin{itemize}

\item {\bf The state space} is $\Delta(\Omega)$; $\mu_t \in \Delta(\Omega)$ is often called the {\it belief}-state. We endow this space with the weak convergence topology, and we note that since $\Omega$ is finite, the set of probability measures on $\Omega$ is a compact space.

\item {\bf The action set} is the set of all maps $\Gamma^1:= \{\gamma^1: \Omega \to \mathbb{A}^1\}.$ We note that since the commitment type policies are given a priori, one could also regard the action set to be the set $\mathbb{A}^1$ itself.\footnote{We note that randomized strategies may also be considered by adding a randomization variable.}

\item {\bf The transition kernel} is given by $P: \Delta(\Omega)\times \Gamma^1 \rightarrow {\cal B}(\Delta(\Omega))$\footnote{${\cal B}(\Delta(\Omega))$ is the set of all Borel sets on $\Delta(\Omega)$. }  so that for all $B \in {\cal B}(\Delta(\Omega))$: 
\begin{align}\label{walvarRecursion}
&\resizebox{0.98\hsize}{!}{$P\bigg( P_\sigma(\omega | s^2_{[0,t-1]}) \in B \bigg|  P_\sigma(\omega|s^2_{[0,t'-1]}), \gamma^1_{t'}, t' \leq t-1 \bigg)$} \nonumber \\
&\resizebox{0.98\hsize}{!}{$=P\bigg( \bigg\{ { \sum_{a^1_{t-1}} P_\sigma(s^2_{t-1}|a^1_{t-1}) P_\sigma(a^1_{t-1} | \omega, s^2_{[0,t-2]})  P_\sigma(\omega |s^2_{[0,t-2]}) \over  \sum_{a^1_{t-1},\omega} P_\sigma(s^2_{t-1}|a^1_{t-1}) P_\sigma(a^1_{t-1} | \omega, s^2_{[0,t-2]})  P_\sigma(\omega |s^2_{[0,t-2]}) } \bigg\} \in B \nonumber $} \\
& \resizebox{0.98\hsize}{!}{$\quad \quad \quad \quad \quad \quad\quad \quad\quad \quad \quad  \quad \quad  \quad \quad \quad \quad \quad \quad  \bigg|  P_\sigma(\omega | s^2_{[0,t'-1]}), \gamma^1_{t'}, t' \leq t-1 \bigg) \nonumber$} \\
&\resizebox{0.98\hsize}{!}{$=P\bigg( \bigg\{ { \sum_{a^1_{t-1}} P_\sigma(s^2_{t-1}|a^1_{t-1}) P_\sigma(a^1_{t-1} | \omega, s^2_{[0,t-2]})  P_\sigma(\omega |s^2_{[0,t-2]}) \over \sum_{a^1_{t-1},\omega} P_\sigma(s^2_{t-1}|a^1_{t-1}) P_\sigma(a^1_{t-1} | \omega, s^2_{[0,t-2]})  P_\sigma(\omega |s^2_{[0,t-2]}) } \bigg\} \in B \nonumber$} \\
& \qquad \qquad \qquad \qquad \qquad \qquad   \bigg|  P_\sigma(\omega | s^2_{[0,t-2]}), \gamma^1_{t-1} \bigg) 
\end{align}

In the above derivation, we use the fact that the term $P_\sigma(a^1_{t-1}|\omega, s^2_{[0,t-2]})$ is uniquely identified by $P_\sigma(\omega |s^2_{[0,{t-2}]})$ and $\gamma^1_{t-1}$. Here, $\gamma^1_{t-1}$ is the {\it control action}.

\item {\bf The per-stage reward function}, given $\gamma^2_t$, is $U(\mu_t,\gamma^{1}): \Delta(\Omega)\times \Gamma^1 \rightarrow \mathbb{R}$  which is defined as follows
\begin{equation} \resizebox{0.40743256 \textwidth}{!}{$ U(\mu_t,\gamma^{1}) := \sum_{\omega}  P_\sigma(\omega | s^2_{[0,t-1]}) \sum_{\mathbb{A}^1} \bigg( {1}_{\{a^1_t=\gamma^{1}(\omega)\}} u^1(a^1_t,\gamma^2_t(P_\sigma(a^1_t|s^2_{[0,t-1]}),s^2_t)) \bigg)$} \label{CostFn2}
\end{equation}
 where $\mu_t=P_\sigma(\omega |s^2_{[0,t-1]})$. Here, $\gamma^2_t$ is a given measurable function of the posterior $P_\sigma(a^1_t|s^2_{[0,t]})$. We note again that for each Bayesian rational short-lived Player 2 we have
\[\resizebox{0.448 \textwidth}{!}{$\gamma^2_t(P_\sigma(a^1_t|s^2_{[0,t-1]}),s^2_t)) \in \mathrm{arg}\max \bigg(\sum_{a^1} P_\sigma(a^1_t|s^2_{[0,t]}) u^2(a^1,a^2)\bigg)$}. \]
\end{itemize}

Lemma \ref{walrandType0} implies that in the finitely repeated setup, under Assumption \ref{EssUnique}, when Player 2s are Bayesian rational, the long-lived strategic Player 1 can depend his strategy only on Player 2s' posterior belief and time without any loss in payoff performance.

Consider now any Perfect Bayesian Equilibrium where the strategic long-lived Player 1 plays a private strategy, since the strategic long-lived Player 1 cannot have a profitable deviation, the public strategy identified in Lemma \ref{walrandType0} must also give him the same payoff against the given sequence of strategies of Player 2s. Hence, in the finitely repeated setup, under Assumption \ref{EssUnique}, any {\bf Perfect Bayesian Equilibrium} payoff of the normal type Player 1, is also a {\bf Perfect Public Equilibrium} payoff.\footnote{A Perfect Public Equilibrium is a Perfect Bayesian Equilibrium where each player uses a public strategy, i.e., a strategy that only depends on the information which is available to both players.} Therefore, given our Markov chain construction: 

\begin{theorem} \label{pbempef}
	In the finitely repeated game, under Assumption \ref{EssUnique}, the set of Perfect Bayesian Equilibrium payoffs of the strategic long-lived Player 1 is equal to the set of Markov Perfect Equilibrium payoffs.
\end{theorem}

\noindent{\bf Proof.} Markov Perfect Equilibrium payoff set is a subset of Perfect Bayesian Equilibrium payoff set. Hence, it is enough to show that for each Perfect Bayesian Equilibrium there exists a properly defined Markov Perfect Equilibrium which is payoff equivalent for the strategic long-lived Player 1. This follows from Lemma \ref{walrandType0} and our Markov chain construction.\qed \\

Lemma \ref{witRealCoding} and Lemma \ref{walrandType0} above have a coding theoretic flavor: The classic works by Witsenhausen \cite{Witsenhausen} and  Walrand and Varaiya \cite{WalrandVaraiya}, are of particular relevance; Teneketzis \cite{Teneketzis} extended these approaches to the more general setting of non-feedback communication and \cite{YukIT2010arXiv} and \cite{YukselBasarBook} extended these results to more general state
spaces (including $\mathbb{R}^d$). Extensions to infinite horizon stages have been studied in \cite{YukLinZeroDelay}. In particular, Lemma \ref{witRealCoding} can be viewed as a generalization of Witsenhausen \cite{Witsenhausen}. On the other hand, Lemma \ref{walrandType0} can be viewed as a generalization of Walrand and Varaiya \cite{WalrandVaraiya} and \cite{YukLinZeroDelay}. The proofs build on \cite{YukIT2010arXiv}. However, these results are different from the above contributions due to the fact that the utility functions do not depend explicitly on the type of Player 1, but depend explicitly on the actions $a^1_t$ and that these actions are not available to Player 2 unlike the setup in  \cite{YukIT2010arXiv}. Next, we consider the infinitely repeated setup in the following.

\subsection{Infinite Horizon and Equilibrium Strategies}\label{InfiniteHorizon}

We proceed with Lemma \ref{walrandType} which is the extension of Lemma \ref{walrandType0} to the infinitely repeated setup. Lemma \ref{walrandType} will be the key result that gives us a similar controlled Markov chain construction for the infinitely repeated game, hence a payoff-equivalent stationary Markov Perfect Equilibrium for each Perfect Bayesian Equilibrium. 

\begin{lemma}\label{walrandType}
In the infinitely repeated game, under Assumption \ref{EssUnique}, given any arbitrary sequence of strategies of Bayesian rational short-lived Player 2s, for any (private) strategy of the strategic long-lived  Player 1, there exists a (public) strategy that only conditions on  $P_\sigma(\omega | s^2_{[0,t-1]}) \in \Delta(\Omega)$ and $t$  which yields the strategic long-lived Player 1 a better payoff against the given sequence of strategies of Player 2s.

Furthermore, the strategic long-lived Player 1's optimal stationary strategy against this given sequence of strategies of Player 2s can be characterized by solving an infinite horizon discounted dynamic programming problem.
\end{lemma}

\noindent{\bf Proof.} See the Appendix. \\

Therefore,  in the infinitely repeated setup as well, under Assumption \ref{EssUnique}, any private strategy of the normal type Player 1 can be replaced, without any loss in payoff performance, with a public strategy which only depends on $P_\sigma(\omega | s^2_{[0,t-1]})$ and $t$. Hence, for any {\bf Perfect Bayesian Equilibrium} there exists a {\bf Perfect Public Equilibrium} which is payoff-equivalent for the strategic long-lived Player 1 in the infinitely repeated game as well.

Furthermore, since there is a stationary optimal public strategy for the strategic long-lived Player 1 against any given sequence of strategies of Bayesian rational Player 2s, any payoff the strategic long-lived Player 1 obtains in a {\bf Perfect Bayesian Equilibrium}, he can also obtain in a {\bf Markov Perfect Equilibrium}.\footnote{A Markov Perfect Equilibrium is a Perfect Bayesian equilibrium where there is a payoff-relevant state space and both players are playing Markov strategies that only depend on the state variable.} 

\begin{theorem} \label{pbempe}
In the infinitely repeated game, under Assumption \ref{EssUnique}, the set of Perfect Bayesian Equilibrium payoffs of the strategic long-lived Player 1 is equal to the set of Markov Perfect Equilibrium payoffs.
\end{theorem}

\noindent{\bf Proof.} The proof follows from Lemma \ref{walrandType} and our Markov chain construction as in the proof of Theorem \ref{pbempef}.\qed \\

Observe that $\{\mu_t(\bar{\omega}) = \mathbb{E}[1_{\omega = \bar{\omega}} | s^2_{[0,t]}]\}$, for every fixed $\bar{\omega}$, is a bounded martingale sequence adapted to the information at Player 2, and as a result as $t \to \infty$, by the submartingale convergence theorem \cite{BorkarBook} there exists (a random) $\bar{\mu}$ such that $\mu_t \to \bar{\mu}$ almost surely. Let $\bar{\mu}$ be an {\bf invariant posterior}, that is, a (sample-path) limit of the $\mu_t$ process. Equation (\ref{DP}) leads to the following fixed point equation:\footnote{Equation (\ref{DP}) appears in the proof of Lemma \ref{walrandType} in the Appendix.}
\[V^1(\omega,\bar{\mu}) = \max_{a^1 = \gamma^1_t(\mu,\omega)}( \mathbb{E}[u^1(a^1_t,\gamma^2(\mu)) + \delta \mathbb{E}[V^1[(\omega,\bar{\mu})]  )\]
Therefore,
\begin{eqnarray}\label{invDisc}
V^1(\omega,\bar{\mu}) = {1 \over 1 - \delta} \max_{\gamma^1_t} \mathbb{E}[u^1(a^1_t,a^2_t(\bar{\mu}))],
\end{eqnarray}
and since the solution is asymptotically stationary, the optimal strategy of the strategic long-lived Player 1 when $\mu_0=\bar{\mu}$ has to be a Stackelberg solution for a Bayesian game with prior $\bar{\mu}$; thus, {\it a Perfect Bayesian Equilibrium strategy for the strategic long-lived Player 1 has to be mimicking the stage game Stackelberg type forever}. This insight will be useful in the following section with further refinements.

\section{Undiscounted Average Payoff Case and An  Upper Payoff Bound for the Arbitrarily Patient Long-lived Player}\label{AverageSection}

We next analyze the setup where the strategic long-lived Player 1 were to maximize his {\bf undiscounted} average payoff instead of his discounted average payoff. Not only we identify an optimal strategy for the strategic long-lived Player 1 in this setup, but also we establish an {\bf upper payoff bound} for the arbitrarily patient strategic long-lived Player 1 in the standard {\bf discounted} average payoff case\textemdash through an Abelian inequality.\footnote{Even though there is a large literature on repeated games with incomplete information in the undiscounted setup, the only papers that we know of that study the reputation games explicitly in the this setup are \cite{CrippsThomas95} and \cite{sorin99}. As opposed to our model, \cite{CrippsThomas95} analyzes a two-person reputation game where both of the players are long-lived. On the other hand, \cite{sorin99} unifies results from merging of probabilities, reputation, and repeated games with incomplete information in both discounted and undiscounted setups.} 

The only difference from our original setup is that the strategic long-lived Player 1 now wishes to maximize
\[ \liminf_{N \to \infty} {1 \over N} \mathbb{E}_{\sigma^1, \sigma^2}^{\mu_0} [\sum_{t=0}^{N-1} u^1(a^1_t,a^2_t)]. \]
Therefore, in any Perfect Bayesian Equilibrium, same as before,  the short-lived (Bayesian rational) Player 2s will continue to be best replying to their updated beliefs. On the other hand, the strategic long-lived Player 1 will be playing a strategy which maximizes his undiscounted average payoff given that each Player 2 will be best replying to their updated beliefs.

The main problem in analyzing the undiscounted setup is that most of the structural coding/signaling results that we have for finite horizon or infinite horizon discounted optimal control problems do not generalize for the undiscounted case, since the construction of controlled Markov chains (which is almost given apriori in stochastic control problems) is based on backwards induction arguments leading to structural results that are applicable only for finite horizon problems. %, see \cite{YukLinZeroDelay}. Therefore, we will arrive at the following results using an indirect approach which is based on more intricate arguments.

Let us re-visit the discounted setup: Let $\bar{\mu}$ be an {\bf invariant posterior}, that is, a (sample-path) limit of the $\mu_t$ process which exists by the discussion with regard to the submartingale convergence theorem. Equation (\ref{invDisc}) is applicable for every $\delta \in (0,1)$ so that
\begin{eqnarray} \label{invDisc2}
(1-\delta)V^1(\omega,\bar{\mu}) = \max_{\gamma^1_t} \mathbb{E}[u^1(a^1_t,a^2_t(\bar{\mu}))],
\end{eqnarray}
and the optimal strategy of the strategic long-lived Player 1 when $\mu_0=\bar{\mu}$ is a Stackelberg solution for a Bayesian game with prior $\bar{\mu}$; thus, {\it a Perfect Bayesian Equilibrium strategy for the strategic long-lived Player 1 has to be mimicking the stage game Stackelberg type forever}. In the following, we will identify conditions when the limit $\bar{\mu}$ will turn out to be a dirac delta distribution at the normal type, that is $\bar{\mu} = \delta_{w}$ (basically, as in the complete information case). 
%In particular, every optimal strategy should be such that if Player 2's belief has converged, then the equilibrium behaviour must be of stage-wise Stackelberg type. Note also that, by the analysis in the previous section, Player 2 behaves as if his strategy is optimal once his opinion is within a neighbourhood of the limit belief. Once Player 2s start best replying to such a neighborhood of this limit belief, Player 1's optimal strategy becomes the Stackelberg action which is maximized according to the limit belief of Player 2s. 
Furthermore, the above discussion implies the following observation: By a direct application of the Abelian inequality (see \ref{Tauberian}), we have that when $\mu_0 = \bar{\mu}$,
\begin{eqnarray}\label{whenInvariant}
&& \sup_{\sigma^1,\sigma^2} \liminf_{N \to \infty} {1 \over N} \mathbb{E}_{\sigma^1,\sigma^2} \bigg[ \sum_{m=0}^{N-1} u^1(a^1_m,a^2_m)\bigg] \nonumber \\
&& \leq \limsup_{\delta \to 1} \sup_{\sigma^1,\sigma^2}  \mathbb{E}_{\sigma^1,\sigma^2} (1 - \delta) \bigg[\sum_{m=0}^{\infty} \delta^m u^1(a^1_m,a^2_m)\bigg]  \nonumber \\
&& = \max_{\gamma^1_t} \mathbb{E}[u^1(a^1_t,a^2_t(\bar{\mu}))],
\end{eqnarray}
where the last equality follows from (\ref{invDisc2}). In the following, we will elaborate further on these observations and arrive at more refined results. We state the following identifiability assumption.
\begin{assumption}\label{unifIdentifiability}
Uniformly over all stationary and optimal (for sufficiently large discount parameters $\delta$) strategies $\tilde{\sigma}^1,\tilde{\sigma}^2$,
\begin{align}
& \lim_{\delta \to 1} \sup_{\tilde{\sigma}^1,\tilde{\sigma}^2} \bigg| \mathbb{E}_{\tilde{\sigma}^1,\tilde{\sigma}^2}(1 - \delta) \bigg[\sum_{t=0}^{\infty} \delta^t u^1(a^1_t,a^2_t)\bigg] \nonumber \\
& \quad \quad -  \limsup_{N \to \infty} {1 \over N}  \mathbb{E}_{\tilde{\sigma}^1,\tilde{\sigma}^2} \bigg[ \sum_{t=0}^{N-1} u^1(a^1_t,a^2_t)\bigg] \bigg| = 0 \label{unifS}
\end{align}
\end{assumption}

A sufficient condition for Assumption \ref{unifIdentifiability} is the following.

\begin{assumption}\label{unifIdentifiabilityI}
Whenever the strategic long-lived Player 1 adopts a stationary strategy, for any initial commitment prior, there exists a stopping time $\tau$ such that for $t \geq \tau$, Player 2s' posterior beliefs become so that his best response does not change (that is, his best-response to his beliefs leads to a constant action). Furthermore, $\mathbb{E}[\tau] < \infty$, uniformly over any stationary strategy $\sigma^1$.
\end{assumption}

%{\it{\color{red} OLD:
%Furthermore, Proposition \ref{IdenCon} below shows that Assumption \ref{unifIdentifiability} is indeed implied by one of the most standard identifiability assumptions in the repeated games literature:
%
%\begin{proposition}\label{IdenCon}
%Consider the matrix $A$ whose rows consist of the vectors:
%\[\resizebox{0.5\textwidth}{!}{$\begin{bmatrix} P(s^2_t=k | a^1_t=1) & P(s^2_t=k | a^1_t=2) & \cdots & P(s^2_t=k | a^1_t=|\mathbb{A}^1|) \end{bmatrix}$} \]
%where $ k \in \{1,2,\cdots, |\mathbb{S}^2|\}$
%If $rank(A)= |\mathbb{A}^1|$, then Assumption \ref{unifIdentifiability} holds.
%\end{proposition}
%}}

Furthermore, Proposition \ref{IdenCon} below shows that Assumption \ref{unifIdentifiability} is indeed implied by one of the most standard identifiability assumptions in the repeated games literature:
%We will also utilize the following {\it observability} condition:
\begin{assumption}\label{observabilityC}
Consider the matrix $A$ whose rows consist of the vectors:
\[\resizebox{0.5\textwidth}{!}{$\begin{bmatrix} P_\sigma(s^2_t=k | a^1_t=1) & P_\sigma(s^2_t=k | a^1_t=2) & \cdots & P_\sigma(s^2_t=k | a^1_t=|\mathbb{A}^1|) \end{bmatrix}$} \]
where $ k \in \{1,2,\cdots, |\mathbb{S}^2|\}$. We have that $rank(A)= |\mathbb{A}^1|$
\end{assumption}

\begin{proposition}\label{IdenCon}
Under Assumption \ref{observabilityC},
\[\|P_\sigma(a^1_t \in \cdot | h^2_t) - P_\sigma(a^1_t \in \cdot | h^2_t, \omega)\|_{TV} \to 0.\]
for every $\sigma$. 
Furthermore, under Assumption \ref{observabilityC}, Assumption \ref{unifIdentifiability} holds.
\end{proposition}

\noindent{\bf Proof.}   See the Appendix. \\

The sufficient condition described in Proposition \ref{IdenCon} is a standard identifiability assumption, sometimes referred as the full-rank monitoring assumption in the reputations literature, see for example \cite[Assumption 2]{crippsmailathsamuelson04}. Under Assumption \ref{unifIdentifiability}, we establish that mimicking a Stackelberg commitment type forever is an optimal strategy for the strategic long-lived Player 1 in the undiscounted setup:

%We also state the following as a separate assumption. Note that for discounted cost problems, we were able to arrive at this assumption without any loss in performance.
%\begin{assumption}\label{MarkovStrategy}
%The strategy of the long-term player is of the type given in Lemma \ref{walrandType0}.
%\end{assumption}

\begin{theorem}\label{ACOE}
In the undiscounted setup, under Assumption \ref{observabilityC}, an optimal strategy for the strategic long-lived Player 1 in the infinitely repeated game is the stationary strategy {\it mimicking the Stackelberg commitment type forever.} 
\end{theorem}
\noindent{\bf Proof.}   See the Appendix. \\
\begin{remark}
\begin{itemize}
\item[(i)] We note that we cannot directly use the arguments in \cite{YukLinZeroDelay} with regard to the optimality of Markovian strategies (those given in Lemma \ref{walrandType0}) for average-cost/average-payoff problems since a crucial argument in that paper is to establish a nearly optimal coding scheme which uses the fact that more information cannot hurt both the encoder and the decoder; in our case here, we have a game and the value of (or the lack of) information can be positive or negative in the absence of a further analysis. 
\item[(ii)] Under the conditions noted, it follows that Player 1 cannot {\it abuse} his reputation in the undiscounted setup: An optimal policy is an {\it honest} stage-wise Stackelberg policy. Abusing (through {\it exploiting}) the reputation is inherently a discounted optimality phenomenon.
%\item[(iii)] It would be natural to expect that some of the conditions we have imposed for the average cost problem can be eliminated. It would be desirable in particular to eliminate the continuity condition; continuity around the invariant belief could be sufficient and a natural condition
\end{itemize}
\end{remark}
As an implication of Theorem \ref{ACOE}, we next state the aforementioned upper bound for Perfect Bayesian Equilibrium payoffs of the arbitrarily patient strategic long-lived Player 1 in the discounted setup as Theorem \ref{upperbound}.

\begin{theorem}\label{upperbound}
Under the assumptions of Theorem \ref{ACOE}, 
\[ \limsup_{\delta \to 1} V^1_{\delta}(\omega,\mu^0) \leq \max_{\alpha_1\in \Delta(A_1), \alpha_2 \in BR(\alpha_1)}u_1(\alpha_1,\alpha_2).\] That is, an upperbound for the value of the reputation for an arbitrarily patient strategic long-lived Player 1 in any Perfect Bayesian Equilibrium of the discounted setup is his stage game Stackelberg equilibrium payoff.
\end{theorem}

%\textbf{Proof of Theorem \ref{upperbound}.} Note the following Abelian inequalities: Let $a_n$ be a sequence of non-negative numbers and $\beta \in (0,1)$. Then,
%\begin{eqnarray}\label{Tauberian} \index{TauberienTheorem}
%&&\liminf_{N \to \infty} {1 \over N} \sum_{m=0}^{N-1} a_m  \leq \liminf_{\beta \uparrow 1} (1 - \beta) \sum_{m=0}^{\infty} \beta^m a_m \nonumber \\
%&&\leq \limsup_{\beta \uparrow 1} (1 - \beta) \sum_{m=0}^{\infty} \beta^m a_m  \leq \limsup_{N \to \infty} {1 \over N} \sum_{m=0}^{N-1} a_m
%\end{eqnarray}
%Therefore, for any $\delta$, an upper bound is obtained by the corresponding undiscounted average payoff problem. Since for every $\delta$, an optimal strategy is stationary, and under the stationary strategy the average payoff converges to the one achieved by the case where the type of Player 1 is correctly identified by Player 2s, the result follows from Theorem \ref{ACOE}.
%\qed

Theorem \ref{upperbound} provides an upper bound on the value of reputation for the strategic long-lived Player 1 in the discounted setup. That is, in the discounted setup, an arbitrarily patient strategic long-lived Player 1 cannot do any better than his best Stackelberg payoff under reputational concerns as well. This upperbound coincides with those provided before by Fudenberg and Levine \cite{fudenberglevine92} and Gossner \cite{gossner11}.

\section{Continuity of Payoff Values}\label{contValuePrior}
Next, we consider the continuity of the payoff values of the strategic long-lived Player 1 in the prior beliefs of Player 2s for any Markov Perfect Equilibrium obtained through the aforementioned dynamic programming. In this section, we assume the following.

\begin{assumption}\label{AssCont}
Either Assumption \ref{EssUnique}(i) holds or the stage game payoff functions are identical for both players. 
\end{assumption}

\begin{lemma}\label{weakContKernel}
The transition kernel of the aforementioned Markov chain is weakly continuous in the (belief) state and action.
\end{lemma}
\noindent{\bf Proof.} See the Appendix. \\

We note that, as in \cite{YukLinZeroDelay} if the game is an identical interest game, the continuity results would follow. By Assumption \ref{AssCont}, the per-stage reward function, $U(\mu_t,\gamma^{1})$, is continuous in $\mu_t$. The continuity of the transition kernel and per-stage reward function together with the compactness of the action space leads to the following continuity result.

\begin{theorem}\label{contOptimalPolicy}
Under Assumption \ref{AssCont}, the value function $V^1_t$ of the dynamic program given  in (\ref{DP}) is continuous in $\mu_t$ for all $t\geq 0$.\footnote{The dynamic program (\ref{DP}) appears in the proof of Lemma \ref{walrandType} in the Appendix.}
\end{theorem}

\textbf{Proof of Theorem \ref{contOptimalPolicy}.}
Given Lemma \ref{weakContKernel} and Assumption \ref{EssUnique}(i), the proof follows from an inductive argument and the measurable selection hypothesis. In this case, the discounted optimality operator becomes a contraction mapping from the Banach space of continuous functions on $\Delta(\Omega)$ to itself, leading to a fixed point in this space. 
\qed\\

Theorem \ref{contOptimalPolicy} implies that any Markov Perfect Equilibrium payoff of the strategic long-lived Player 1 obtained through the dynamic program in (\ref{DP}) is robust to small perturbations in the prior beliefs of Player 2s under Assumption \ref{EssUnique}. This further implies that the following conjecture made by Cripps, Mailath, and Samuelson \cite{crippsmailathsamuelson04} is indeed true in our setup: There exists a particular equilibrium in the complete information game and a bound such that for \emph{any} commitment type prior (of Player 2s) less than this bound, there exists an equilibrium of the incomplete information game where the strategic long-lived Player 1's payoff is arbitrarily close to his payoff from the particular equilibrium in the complete information game.\footnote{This conjecture appears as a presumption of \cite[Theorem 3]{crippsmailathsamuelson04}. They write ``We conjecture this hypothesis is redundant, given the other conditions of the theorem, but have not been able to prove it''.} This is also in line with the findings of \cite{Dalkiran16}, which uses the methods of \cite{abpest90} to show a similar upper semi continuity result.

For the undiscounted setup, however, in Section \ref{AverageSection}, we were able to achieve a much stronger continuity result, without requiring Assumption \ref{AssCont} but instead Assumption \ref{observabilityC}, in addition to the assumptions stated at the beginning of the paper. We formally state this result next.
\begin{theorem}\label{contOptimalPolicyAverage}
Under the conditions of Theorem \ref{ACOE}, the undiscounted average value function does not depend on the prior $\mu_0$.
\end{theorem}

\section{A Lower Payoff Bound on Reputation through Measure Concentration}\label{lowerBSection}

We next identify a lower payoff bound for the value of reputation through an explicit measure concentration analysis. As mentioned before, it was Fudenberg and Levine \cite{fudenberglevine89}, \cite{fudenberglevine92} who provided such a lower payoff bound for the first time, to our knowledge. They constructed a lower bound for any equilibrium payoff of the strategic long-lived player by showing that Bayesian rational  short-lived players can be surprised at most finitely many times when a strategic long-lived player mimics a commitment type forever. Using the chain rule property of the concept of relative entropy, \cite{gossner11} obtained a lower bound for any equilibrium payoff of the strategic long-lived player by showing that any equilibrium payoff of the strategic long-lived player is bounded from below (and above) by a function of the average discounted divergence between the prediction of the short-lived players conditional on the long-lived player's type and its marginal. 

Our analysis below provides a sharper lower payoff bound for the value of reputation through a refined measure concentration analysis. To obtain this lower bound, as in \cite{fudenberglevine92} as well as \cite{gossner11}, we let the strategic long-lived Player 1 mimic (forever) a commitment type, $\hat{\omega}=m$, to investigate the best responses of the short-lived Player 2s. In any Perfect Bayesian Equilibrium, such a deviation, i.e., deviating to mimicking a particular commitment type forever, is always possible for the strategic long-lived Player 1.

Let $|\Omega|=M$ be the number of all possible types of the long-lived Player 1. {We will assume for simplicity that all the types are deterministic, as opposed to the more general mixed types considered earlier in the paper}. With $m$ being the type mimicked forever by Player 1, we will identify  a function $f$ below such that for any $\hat{\omega} \in \hat{\Omega}$ when criterion (\ref{criterion}) below holds,
\begin{eqnarray}\label{criterion}
{P_\sigma(\omega = m | s^2_{[0,t]}) \over P_\sigma(\omega=\hat{\omega} | s^2_{[0,t]}) } \geq f(M),
\end{eqnarray}
Player 2 of time $t$ will act as if he knew the type of the long-lived Player 1 is $m$. This will follow from the fact that $\max_{a^2} \sum P_\sigma(\hat{\omega}| s^2_{[0,t]}) u^2(a^1,a^2)$ is continuous in $P_\sigma(\hat{\omega}| s^2_{[0,t]})$ and that $P_\sigma(\hat{\omega}| s^2_{[0,t]})$ concentrates around the true
type under a mild informativeness condition on the observable variables.
%Let
%\begin{eqnarray}
%&&\tau_{m} = \min\{T \geq 0: \max_{a^2} \sum_{a^1} P_\sigma(a^1| s^2_{[0,t]}) u^2(a^1,a^2)  \nonumber \\
%&& \quad \quad \quad \quad = \max_{a^2} \sum_{a^1} P_\sigma(a^1| \omega=m) u^2(a^1,a^2) \quad \forall  t \geq T\}, \nonumber
%\end{eqnarray}
%Intuitively, $\tau_{m}$ is a time when Player 2 behaves as if the type of the long-lived Player 1 is $m$ as far as their optimal strategies are concerned. 
Let
\begin{eqnarray}
&&\tau_{m} = \{ t \geq 0: \max_{a^2} \sum_{a^1} P_\sigma(a^1| s^2_{[0,t]}) u^2(a^1,a^2)  \nonumber \\
&& \quad \quad \quad \quad = \max_{a^2} \sum_{a^1} P_\sigma(a^1| \omega=m) u^2(a^1,a^2)\} \nonumber
\end{eqnarray}
Intuitively, $\tau_{m}$ is the (random) set of times that Players 2 behave as if the type of the long-lived Player 1 is $m$ as far as their optimal strategies are concerned. 

\begin{lemma}\label{boundLearning}
Let $\epsilon > 0$ be such that for any $\bar{a}^1 \in \mathbb{A}^1$ and $\tilde{a}^2, \hat{a}^2 \in \mathbb{A}^2$
\[|u^2(\bar{a}^1,\tilde{a}^2) - u^2(\bar{a}^1, \hat{a}^2)|\geq {\epsilon \over 1 - \epsilon} \bigg(\max_{a^1,a^2} |u^2(a^1,a^2)| \bigg)\]
If (\ref{criterion}) holds at time $t$  when  $f(M)= { (1 - \epsilon) \over \epsilon}M$, then $t \in \tau_{m}$.
\end{lemma}

\noindent{\bf Proof.} See the Appendix. \\

Lemma \ref{boundLearning} implies that when criterion  (\ref{criterion}) holds to be true for $f(M)= { (1 - \epsilon) \over \epsilon}M$, at time $t$ any Player 2 of time $t$ and onwards will be best responding to the commitment type $m$. This can be interpreted as the long-lived Player {\bf having a reputation to behave like type $m$} when criterion (\ref{criterion}) is satisfied.

\begin{theorem}\label{geometricBound}
Suppose that $0 < {P_\sigma(s^2|\omega=m) \over P_\sigma(s^2|\omega=\hat{\omega})} < \infty$ for all $\hat{\omega} \in \hat{\Omega}$ and $s^2 \in \mathbb{S}^2$. For all $k \in \mathbb{N}$, $P_\sigma(k \notin \tau_{m}) \leq R \rho^{k}$ for some $\rho \in (0,1)$ and $R \in \mathbb{R}$.
\end{theorem}

\noindent{\bf Proof.} See the Appendix. \\

We are now ready to provide our lower  bound for Perfect Bayesian Equilibrium payoffs of the strategic long-lived Player 1, for a fixed discount factor $\delta \in(0,1)$.

\begin{theorem}\label{theorembound}
A lower bound for the expected payoff of the strategic long-lived Player 1 in any Perfect Bayesian Equilibrium (in the discounted setup) is given by $\max_{m \in \hat{\Omega}} L(m)$ where
\[\resizebox{0.49\textwidth}{!}{$L(m) = \mathbb{E}_{\{\omega= m\}}\bigg[\sum_{k \notin \tau_{m}} \delta^k u^1(a^1_t,a^2_t)\bigg] + \mathbb{E}_{\{\omega= m\}} \bigg[\sum_{k \in \tau_{m}} \delta^k {\underline{u}^1}^*(m)\bigg]$}\]
where $ {\underline{u}^1}^*(m):=\min_{a^2 \in BR^2(m)} u^1(m,a^2) $ and $BR^2(m):=\mathrm{arg}\max_{a^2 \in \mathbb{A}^2} u^2(m,a^2)$.
\end{theorem}

\noindent{\bf Proof.} By Theorem \ref{geometricBound}, the discounted average payoff can be lower bounded by the sum of the following two terms:
\[ \mathbb{E}_{\{\omega= m\}}\bigg[\sum_{k \notin \tau_{m}} \delta^k u^1(a^1_t,a^2_t)\bigg] + \mathbb{E}_{\{\omega= m\}} \bigg[\sum_{k \in \tau_{m}} \delta^k {\underline{u}^1}^*(m)\bigg]\] where $ {\underline{u}^1}^*(m):=\min_{a^2 \in BR^2(m)} u^1(m,a^2) $ and $BR^2(m):=\mathrm{arg}\max_{a^2 \in \mathbb{A}^2} u^2(m,a^2)$. Since a deviation to mimicking any of the commitment types forever is available to the strategic long-lived Player 1 in any Perfect Bayesian Equilibrium, taking the maximum of the lower bound above for all commitment types gives the desired result. \qed

Observe that when $m$ is a Stackelberg type, i.e., a commitment type who is committed to play the stage game Stackelberg action $
\mathrm{arg}\max_{\alpha_1\in \Delta(A_1)}u_1(\alpha_1,BR^2(\alpha^1))$  for which Player 2s have a unique best reply then \[{\underline{u}^1}^*(m)= \max_{\alpha_1\in \Delta(A_1), \alpha_2 \in BR(\alpha_1)} u_1(\alpha_1,\alpha_2)\]
becomes the stage game Stackelberg payoff.

We next turn to the case of the arbitrarily patient strategic long-lived Player 1. That is, what happens when $\delta \rightarrow 1$. To emphasize the dependence on $\delta$, we use a superscript in $L^{\delta}(m)$.

\begin{theorem} \label{ReputationStackelberg}
\[ \lim_{\delta \to 1} (1 - \delta) L^{\delta}(m) \geq {\underline{u}^1}^*(m) \]
%\sy{TO DO: Verify that the limit exists explicitly.}
\end{theorem}

\noindent\textbf{Proof.} The proof follows from Theorem \ref{theorembound} by taking the limit $\delta \to 1$. Since in $\tau_m$, we can bound the payoff to strategic long-lived Player 1 below by the worst possible payoff, and in $\tau_m$ the strategic long-lived Player 1 guarantees the associated Stackelberg payoff, we obtain by an application of the Abelian inequality, the desired result. 
%
%URBANA
%
% \[\resizebox{1\hsize}{!}{$\lim_{\delta \to 1} L(m) \geq \lim_{\delta \to 1} \mathbb{E}[1-{\delta^{\tau_m}}]  \min_{a^1,a^2} u^1(a^1,a^2) +  \lim_{\delta \to 1} \mathbb{E}[\delta^{\tau_m} {\underline{u}^1}^*(m)] ={\underline{u}^1}^*(m).$}\]
%That $\lim_{\delta \to 1} \mathbb{E}[{\delta^{\tau_m} - 1}]=0$ and $\lim_{\delta \to 1} \mathbb{E}[{\delta^{\tau}}]=1$ follow from the dominated convergence theorem and the fact that $\tau_m$ is finite with probability 1.
 \qed

Theorem \ref{ReputationStackelberg} implies that the lower payoff bound that we provided in Theorem \ref{theorembound} coincides in the limit as $\delta \rightarrow 1$ with those of Fudenberg and Levine \cite{fudenberglevine92} and Gossner \cite{gossner11}. That is, if there exists a Stackelberg commitment type, an arbitrarily patient strategic long-lived Player 1 can guarantee himself a payoff arbitrarily close to the associated Stackelberg payoff in every Perfect Bayesian Equilibrium in the discounted setup.

\section{Conclusion}

In this paper, we studied the reputations problem of an informed long-lived player who controls his reputation against a sequence of uninformed short-lived players by employing tools from stochastic control theory. Our findings contribute to the reputations literature by obtaining new results on the structure of equilibrium behavior in finite-horizon, infinite-horizon, and undiscounted settings, as well as continuity results in the prior probabilities, and improved upper and lower bounds on the value of reputations. In particular, we exhibited that a control theoretic formulation can be utilized to characterize the equilibrium behavior. Even though there are studies that employed dynamic programming methods to study reputation games in the literature, e.g., \cite{jullienpark14},  these studies restrict themselves directly to Markov strategies\textemdash hence to the concept of Markov Perfect Equilibrium without mentioning its relation to the more general (and possibly more appropriate) concept of Perfect Bayesian Equilibrium. Under technical assumptions, we have identified that a nested information structure implies the equivalence of the set of Markov Perfect Equilibrium payoffs and the set of Perfect Bayesian Equilibrium payoffs. It is our hope that the machinery we provide in this paper will open a new avenue for applied work studying reputations in different frameworks.

\appendix
\section{Appendix}

\subsection{Proof of Lemma \ref{witRealCoding}.}
At time $t=T$, the payoff function can be written as follows, where $\gamma^2_t$ denotes a given fixed strategy for Player 2:
\begin{eqnarray}\label{costCompactForm}
\mathbb{E}[u^1(a^1_{t},\gamma^2_{t}(s^2_{[0,t]})) | s^2_{[0,t-1]}]  = \mathbb{E}[F(a^1_t,s^2_{[0,t-1]},s^2_t) | s^2_{[0,t-1]}] \nonumber
\end{eqnarray}
where, $F(a^1_t,s^2_{[0,t-1]},s^2_t) = u^1(a^1_{t},\gamma^2_{t}(s^2_{[0,t]})) $.

Now, by a stochastic realization argument (see Borkar \cite{BorkarRealization}), we can write $s^2_t = R(a^1_t,v_t)$ for some independent noise process $v_t$. As a result, the expected payoff conditioned on $s^2_{[0,t-1]}$ is equal to, by the smoothing property of conditional expectation, the following:
$$\mathbb{E}\bigg[\mathbb{E} [G(a^1_t,s^2_{[0,t-1]},v_t) | \omega, a^1_t,s^2_{[0,t-1]}]  \bigg| s^2_{[0,t-1]}\bigg],$$
for some $G$. Since $v_t$ is independent of all the other variables at times $t' \leq t$, it follows that there exists $H$ so that $\mathbb{E} [G(a^1_t,s^2_{[0,t-1]},v_t) | \omega, a^1_t,s^2_{[0,t-1]}] =: H(\omega, a^1_t,s^2_{[0,t-1]})$.
Note that when $\omega$ is a commitment type, $a^1_t$ is fixed quantity or a fixed random variable.

Now, we will apply Witsenhausen's two stage lemma \cite{Witsenhausen}, to show that we can obtain a lower bound for the double expectation by picking $a^1_t$ as a result of a measurable function of $\omega,s^2_{[0,t-1]}$. Thus, we will find a strategy which only uses $(\omega,s^2_{[0,t-1]})$ which performs as well as one which uses the entire memory available at Player 1. To make this precise, let us fix $\gamma^2_t$ and define for every $k \in \mathbb{A}^1$:
\begin{equation}\label{Witsenconstruct12}
\resizebox{1\hsize}{!}{$\beta_{k} := \bigg\{\omega,s^2_{[0,t-1]} : G(\omega,s^2_{[0,t-1]},k)  \leq G(\omega,s^2_{[0,t-1]},q), \forall q \neq k \bigg\}. \nonumber$}
\end{equation}
 Such a construction covers the domain set consisting of \phantom{X} $(x_t,q_{[0,t-1]})$ but possibly with overlaps. It covers the elements in $\Omega \times \prod_{t=0}^{T-1} \mathbb{S}^2$, since for every element in this product set, there is a maximizing $k \in \mathbb{A}^1$. To avoid the overlap, define a function $\gamma^{*,1}_t$ as:
$$q_t = \gamma^{*,1}_t(\omega,s^2_{[0,t-1]})= k, \quad \rm{if} (\omega,s^2_{[0,t-1]}) \in \beta_{k} \setminus \cup_{i=1}^{k-1} \beta_{i},$$ with $\beta_0=\emptyset$. The new strategy performs at least as well as the original strategy even though it has a restricted structure.

The same discussion applies for earlier time stages as we discuss below. We iteratively proceed to study the other time stages.
For a three-stage problem, the payoff at time $t=2$ can be written as:
\begin{equation}
\resizebox{1\hsize}{!}{$ \mathbb{E}\bigg[ u^1(a^1_2,\gamma^2_2(s^2_1,s^2_2)) + \mathbb{E}[u^1(\gamma^{*,1}_3(\omega,s^2_{[1,2]}),\gamma^2_3\bigg(s^2_1, s^2_2,R(\gamma^{*,1}_3(\omega,s^2_{[1,2]}),v_3)\bigg) | \omega, s^2_1, s^2_2] \bigg| s^2_1 \bigg] \nonumber$}
\end{equation}
The expression inside the expectation is equal to for some measurable $F_2$, $F_2(\omega, a^1_2, s^2_1, s^2_{2})$. Now, once again expressing $s^2_2=R(a^1_2,v_2)$, by a similar argument as above, a strategy at time $2$ which uses $\omega$ and $s^1_2$ and which performs at least as good as the original strategy can be constructed. By similar arguments, a strategy at time $t$, $1 \leq t\leq T$ only uses $(\omega,s^2_{[1,t-1]})$ can be constructed. The strategy at time $t=0$ uses $\omega$.
\hfill \qed

\subsection{Proof of Lemma \ref{walrandType0}.}

The proof follows from a similar argument as that for Lemma \ref{witRealCoding}, except that the information at Player 2 is replaced by the sufficient statistic that Player 2 uses: his posterior information. At time $t=T-1$, an optimal Player 2 will use $P_\sigma(a^1_t|s^2_{[0,t]})$ as a sufficient statistic for an optimal decision. Let us fix a strategy for Player 2 at time t, $\gamma^2_t$ which only uses the posterior $P_\sigma(a^1_t|s^2_{[0,t]})$ as its sufficient statistic.
Let us further note that:
\begin{eqnarray}
&& P_\sigma(a^1_t|s^2_{[0,t]}) = { P_\sigma(s^2_t,a^1_t|s^2_{[0,t-1]}) \over \sum_{a^1_t} P_\sigma(s^2_t,a^1_t|s^2_{[0,t-1]}) } \nonumber \\
&& \resizebox{0.75\hsize}{!}{$= { \sum_{\omega} P_\sigma(s^2_t|a^1_t )P_\sigma(a^1_t | \omega,s^2_{[0,t-1]}) P_\sigma(\omega | s^2_{[0,t-1]})  \over \sum_{\omega} \sum_{a^1_t} P_\sigma(s^2_t|a^1_t )P_\sigma(a^1_t | \omega,s^2_{[0,t-1]}) P_\sigma(\omega | s^2_{[0,t-1]})} \label{whySeparationHolds11}$}
\end{eqnarray}
The term $P_\sigma(a^1_t | \omega, s^2_{[0,t-1]})$ is determined by the strategy of Player 1 (this follows from Lemma \ref{witRealCoding}), $\gamma^1_t$.

As in \cite{YukselBasarBook}, this implies that the payoff at the last stage conditioned on $s^2_{[0,t-1]}$ is given by

\begin{equation}\label{costCompactForm}
\resizebox{1\hsize}{!}{$\mathbb{E}\bigg[u^1\bigg(a^1_{t},\gamma^2_{t}(P_\sigma(a^1_t = \cdot |s^2_{[0,t]}) )\bigg) | s^2_{[0,t-1]}\bigg]  = \mathbb{E}\bigg[F\bigg(a^1_t,\gamma^1_t,P_\sigma(\omega= \cdot | s^2_{[0,t-1]})\bigg) | s^2_{[0,t-1]}\bigg] \nonumber$}
\end{equation}
where, as earlier, we use the fact that $s^2_t$ is conditionally independent of all the other variables at times $t' \leq t$ given $a^1_t$. 
Let $\gamma^{1,s^2_{[0,t-1]}}_t$ denote the strategy of Player 1. The above state is then equivalent to, by the smoothing property of conditional expectation, the following:
\begin{align}
& \resizebox{0.95\hsize}{!}{$ \mathbb{E}\bigg[\mathbb{E} \bigg[ F\bigg(a^1_t,\gamma^1_t,P_\sigma(\omega= \cdot | s^2_{[0,t-1]})\bigg) | \omega, \gamma^{1,s^2_{[0,t-1]}}_t, P_\sigma(\omega= \cdot | s^2_{[0,t-1]}),s^2_{[0,t-1]} \bigg]  \bigg| s^2_{[0,t-1]}\bigg] \nonumber$} \\
& \resizebox{0.95\hsize}{!}{$=  \mathbb{E}\bigg[\mathbb{E} \bigg[ F\bigg(a^1_t,\gamma^1_t,P_\sigma(\omega= \cdot | s^2_{[0,t-1]})\bigg) | \omega, \gamma^{1,s^2_{[0,t-1]}}, P_\sigma(\omega= \cdot | s^2_{[0,t-1]})\bigg]  \bigg| s^2_{[0,t-1]}\bigg] $} \nonumber \\ \label{BayesianBoth}
\end{align}
The second line follows since once one picks the strategy $ \gamma^{1,s^2_{[0,t-1]}}$, the dependence on $s^2_{[0,t-1]}$ is redundant given \\ $P_\sigma(\omega= \cdot | s^2_{[0,t-1]})$.

Now, one can construct an equivalence class among the past $s^2_{[0,t-1]}$ sequences which induce the same $\mu_t(\cdot)= P_\sigma(\omega \in \cdot | s^2_{[0,t-1]})$, and can replace the strategy in this class with one, which induces a higher payoff among the finitely many elements in each class for the final time stage. An optimal output thus may be generated using $\mu_t$ and $\omega$ and $t$, by extending Witsenhausen's argument used earlier in the proof of Lemma \ref{witRealCoding} for the terminal time stage. Since there are only finitely many past sequences and finitely many $\mu_t$, this leads to a (Borel measurable) selection of $\omega$ for every $\mu_t$, leading to a measurable strategy in $\mu_t, \omega$. Hence, the final stage payoff can be expressed as $F_{t}(\mu_t)$ for some $F_t$, without any performance loss.

The same argument applies for all time stages. To show this, we will apply induction as in \cite{YukIT2010arXiv}. At time $t=T-1$, the sufficient statistic both for the immediate payoff, and the {\it continuation payoff} is $P_\sigma(\omega | s^2_{[0,t-1]})$, and thus for the payoff impacting the time stage $t=T$, as a result of the optimality result for $\gamma^1_{T}$. To show that the separation result generalizes to all time stages, it suffices to prove that $\{(\mu_{t},\gamma^1_t)\}$ has a controlled Markov chain form, if the players use the structure above.

Now, for $t \geq 1$, for all $B \in {\cal B}(\Delta(\Omega))$: 

\begin{align}\label{walvarRecursion}
&P\bigg( P_\sigma(\omega | s^2_{[0,t-1]}) \in B \bigg|  P_\sigma(\omega|s^2_{[0,t'-1]}), \gamma^1_{t'}, t' \leq t-1 \bigg) \nonumber \\
&\resizebox{0.98\hsize}{!}{$=P\bigg( \bigg\{ { \sum_{a^1_{t-1}} P_\sigma(s^2_{t-1}|a^1_{t-1}) P_\sigma(a^1_{t-1} | \omega, s^2_{[0,t-2]})  P_\sigma(\omega |s^2_{[0,t-2]}) \over  \sum_{a^1_{t-1},\omega} P_\sigma(s^2_{t-1}|a^1_{t-1}) P_\sigma(a^1_{t-1} | \omega, s^2_{[0,t-2]})  P_\sigma(\omega |s^2_{[0,t-2]}) } \bigg\} \in B \nonumber$} \\
& \resizebox{0.98\hsize}{!}{$\quad \quad \quad \quad \quad \quad\quad \quad\quad \quad \quad  \quad \quad  \quad \quad \quad \quad \quad \quad  \bigg|  P_\sigma(\omega | s^2_{[0,t'-1]}), \gamma^1_{t'}, t' \leq t-1 \bigg) \nonumber$} \\
& \resizebox{0.98\hsize}{!}{$=P\bigg( \bigg\{ { \sum_{a^1_{t-1}} P_\sigma(s^2_{t-1}|a^1_{t-1}) P_\sigma(a^1_{t-1} | \omega, s^2_{[0,t-2]})  P_\sigma(\omega |s^2_{[0,t-2]}) \over \sum_{a^1_{t-1},\omega} P_\sigma(s^2_{t-1}|a^1_{t-1}) P_\sigma(a^1_{t-1} | \omega, s^2_{[0,t-2]})  P_\sigma(\omega |s^2_{[0,t-2]}) } \bigg\} \in B \nonumber$} \\
&  \resizebox{0.98\hsize}{!}{$\quad \quad \quad \quad \quad \quad\quad \quad\quad \quad \quad  \quad \quad  \quad \quad \quad \quad \quad \quad  \bigg|  P_\sigma(\omega | s^2_{[0,t'-1]}), \gamma^1_{t'}, t' = t-1 \bigg)$} \nonumber \\
\end{align}

In the above derivation, we use the fact that the term \\ $P_\sigma(a^1_{t-1}|\omega, s^2_{[0,t-2]})$ is uniquely identified by $P_\sigma(\omega |s^2_{[0,{t-2}]})$ and $\gamma^1_{t-1}$. \qed

\subsection{Proof of Lemma \ref{walrandType}.}
First, going from a finite horizon to an infinite horizon follows from a change of order of limit and infimum as we discuss in the following.
Observe that for any strategy $\{\gamma^1_t\}$ and any $T \in \mathbb{N}$:
\[\mathbb{E}[ \sum_{t=0}^{T-1} \delta^t u^1(a^1_t, a^2_t)] \geq  \inf_{\{\gamma^1_t\}} \mathbb{E}[\sum_{t=0}^{T-1} \delta^t u^1(a^1_t, a^2_t)]\]
and thus
\[ \lim_{T \to \infty} \mathbb{E}[ \sum_{t=0}^{T-1} \delta^t u^1(a^1_t, a^2_t)] \geq \limsup_{T \to \infty} \inf_{\{\gamma^1_t\}} \mathbb{E}[\sum_{t=0}^{T-1} \delta^t u^1(a^1_t, a^2_t)]\]
Since the above holds for an arbitrary strategy, it follows then that
\begin{eqnarray}
&& \inf_{\{\gamma^1_t\}} \lim_{T \to \infty} \mathbb{E}[ \sum_{t=0}^{T-1} \delta^t u^1(a^1_t, a^2_t)] \nonumber \\
&& \geq \limsup_{T \to \infty} \inf_{\{\gamma^1_t\}} \mathbb{E}[\sum_{t=0}^{T-1} \delta^t u^1(a^1_t, a^2_t)] \nonumber 
\end{eqnarray}
On the other hand, due to the discounted nature of the problem, the right hand side can be studied through the dynamic programming (Bellman) iteration algorithms: The following dynamic program holds: Let $\mu_t(w) = P_\sigma(\omega = w| s^2_{[0,t-1]})$.

\begin{eqnarray} \label{DP}
&& V^1(\omega,\mu_t) = \mathbb{T}(V^1)(\omega,\mu_t) := \nonumber \\
&&  \quad \max_{\gamma^1_t} \bigg( \mathbb{E} [u^1(a^1_t,a^2_t) + \delta \mathbb{E}[V^1(\omega,\mu_{t+1}) | \mu_t, \gamma^1_t] \bigg) 
%& \qquad \qquad =:\mathbb{T}(V^1)(\omega,\mu_t)
\end{eqnarray}
where $\mathbb{T}$ is an operator defined by:
\begin{equation}
\resizebox{1\hsize}{!}{$\mathbb{T}(f)(\omega,\mu_t)= \max_{\gamma^1_t}\bigg( \mathbb{E}\bigg[u^1(a^1_t,a^2_t) + \delta \mathbb{E}[f(\omega,\mu_{t+1}) | \mu_t, \gamma^1_t \bigg] \bigg) \nonumber$}
\end{equation}
A value iteration sequence with $V^1_0=0$ and $V_{t+1}=\mathbb{T}(V_t)$, which is well defined by the measurable selection conditions noted in \cite{HernandezLermaMCP} due to the finiteness of our action set (and hence continuity of the iterations in the actions), leads to a stationary solution. This is an infinite horizon discounted payoff optimal dynamic programming equation with finite action spaces (where the strategy is now the {\it action} $\gamma^1_t$). Since the action set is finite in our formulation, it follows that there is a stationary solution as $t \to \infty$. Thus, the sequence of maximizations $\sup_{\gamma^1} \mathbb{E}[\sum_{t=0}^{T-1} \delta^t u^1(a^1_t, a^2_t)]$ leads to a stationary solution as $T \to \infty$, and this sequence of policies admit the structure given in the statement of the theorem. \qed

\subsection{Proof of Proposition \ref{IdenCon}.} 
{Recall the the chain rule of relative entropy implies the following: For joint measures $P,Q$ on random variables $X,Y$ with finite relative entropy, we have $D(P(X,Y)\|Q(X,Y))=D(P(X)\|Q(X))+D(P(Y|X)\|Q(Y|X))$. Let $X=\omega$ and $Y:=s^2_{[0,\infty)}$, $P:=P_{\sigma,\omega=w}((\omega,s^2_{[0,\infty)}) \in \cdot)$ (i.e., with the true distribution given the type of the long-run player) and $Q:=P_\sigma((\omega,s^2_{[0,\infty)}) \in \cdot)$ (this is the distribution seen by Players 2). Then  (following \cite{gossner11}, see also \cite[Section 8]{McDonaldYuksel}) the conditional relative entropies are summable with the bound $D({\delta_{\omega}} | \mu_0) < \infty$, which also implies that
\[\mathbb{E}\bigg[D\bigg( P_\sigma(s^2_t \in \cdot | h^2_t,\omega) \bigg|\bigg| P_\sigma(s^2_t \in \cdot | h^2_t) \bigg)\bigg] \to 0.\]
%\sy{ELABORATE ON THIS FURTHER; THE SETUP IS SLIGHTLY DIFFERENT FROM GOSSNER'S SETTING, PERHAPS WRITE THE CHAIN RULE EXPLICITLY. PERHAPS UNIFORMITY WILL COME DIRECTLY.}
From Pinsker's inequality noting that convergence in total variation is implied by convergence in relative entropy:
\begin{eqnarray}\label{conver}
\mathbb{E}[||P_\sigma(s^2_t \in \cdot | h^2_t) - P_\sigma(s^2_t \in \cdot | h^2_t, \omega)||_{TV}^2] \to 0
\end{eqnarray}
where the expectation is with respect to the true distribution (given the type of the long-run player). But,
\[P_\sigma(s^2_t =s | h^2_t) = \sum_{a^1} P_\sigma(s^2_t=s | a^1_t=a^1) P_\sigma(a^1_t=a^1 | h^2_t)\]
Thus, all we need to ensure is that Player 2's belief $P_\sigma(a^1_t \in \cdot | h^2_t)$ is sufficiently close to a terminal value. Suppose that the conditions of the theorem holds, but $|P_\sigma(a^1_t | h^2_t) - P_\sigma(a^1_t | h^2_t, \omega)|$ $>$ $\delta$ for some subsequence of time values. If the rank of $A$ is $|\mathbb{A}^1|$, then, $|P_\sigma(a^1_t | h^2_t) - P_\sigma(a^1_t | h^2_t, \omega)| > \delta$ would imply that $|P_\sigma(s^2_t | h^2_t) - P_\sigma(s^2_t | h^2_t, \omega)| > \epsilon$ for some positive $\epsilon$, which would be a contradiction (to see this, observe that the vector $P_\sigma(a^1_t \in \cdot  | h^2_t) - P_\sigma(a^1_t \in \cdot  | h^2_t, \omega)$ cannot be orthogonal to each of the rows of $A$, due to the rank condition). In particular, (\ref{conver}) implies the convergence of $P_\sigma(a^1_t \in \cdot | h^2_t) - P_\sigma(a^1_t \in \cdot | h^2_t, \omega)$ to zero: the summability (bounded from above) of the conditional relative entropies implies that the expected number of instances where the error between the conditional probabilities is above any specified amount will be finite (uniform over all policies).} 

%{\it{\color{red} PROOF FOR: Furthermore, under Assumption \ref{AssCont}, Assumption %\ref{unifIdentifiability} holds.}
Now, using uniform continuity of the per-stage utility in the posterior of player 2 (e.g., through a related result from Gossner \cite{gossner11}), we can uniformly bound the error in the per-stage from the setup when the posterior seen by Player 2 is exactly $P_\sigma(a^1_t \in \cdot  | h^2_t, \omega)$ (where crucially the error is uniform over all posteriors, regardless of the strategy of Player 1). In particular, the pay-off into the future would be so that, it would be within the pay-off for the setup when the posterior of player 2 would correspond to having the prior $\delta_{w}$ on the normal type, as in the complete information case, for any considered normal policy (which in the statement of Assumption \ref{unifIdentifiability} is a stationary policy). On the other hand, we know, by the analysis in (\ref{invDisc2}) that any optimal stationary policy with prior $\delta_w$ will be a stage-wise Stackelberg policy, and the average pay-off (the right-hand side of (\ref{unifS}) in this case will correspond exactly to \[(1-\delta)V^1(\omega,\delta_w) = \max_{\gamma^1_t} \mathbb{E}[u^1(a^1_t,a^2_t(\delta_w))].\]
Together with the uniformity (over strategies) of the relative entropy bound, we conclude that Assumption \ref{unifIdentifiability} holds.

%In particular, the type of Player 1 will be essentially known by Players 2. Now, if the type of Player 1 is known, then the discounted optimal policies will be of commitment types maximizing the stage-wise Stackelberg policy, which in turn will be mimicking a Stackelberg commitment type. This then implies Assumption \ref{unifIdentifiability}.} %CHECK: This would imply that Assumption \ref{unifIdentifiability} would hold if all types were commitment, we need to think if we wish to claim the uniformity over any stationary normal type, however). }
% \sy{Furthermore, the summability due to the chain rule of relative entropy implies that the relative entropy will be above any fixed quantity only finitely many times in expectation, uniform over all policies. SHALL WE ADD THIS? THIS WILL INDEED IMPLY UNIFORMITY TOO, AS ORIGINALLY CLAIMED.}
%\sy{PLAYERS 2 CAN BE OUTSIDE THE NEIGHBORHOOD OF THE INFINITE LIMIT ONLY FINITELY MANY TIMES, SINCE THE RELATIVE ENTROPY CAN BE DIFFERENT ONLY FINITELY MANY TIMES DUE TO THE SUMMABILITY CONDITION. WE CAN TAKE OUT THE REST I THINK: Furthermore, Jensen's inequality implies that
%\begin{eqnarray}\label{conver2}
%\resizebox{0.93\hsize}{!}{$|\mathbb{E}[P_\sigma(s^2_t =s | h^2_t) - P_\sigma(s^2_t =s | h^2_t, \omega)]| \leq  \mathbb{E}[|P_\sigma(s^2_t =s | h^2_t) - P_\sigma(s^2_t =s | h^2_t, \omega)|] \to 0$}
%\end{eqnarray}
%and thus in finite expected time the deviation in the conditional probabilities will be less than a prescribed amount and Assumption \ref{unifIdentifiability} holds.}
\qed

\subsection{Proof of Theorem \ref{ACOE}}

Note the following {\it Abelian} inequalities (see, e.g., Lemma 5.3.1 in Hernandez-Lerma and Lasserre \cite{HernandezLermaMCP}): Let $a_n$ be a sequence of non-negative numbers and $\beta \in (0,1)$. Then,
\begin{align}\label{Tauberian} \index{TauberienTheorem}
&\liminf_{N \to \infty} {1 \over N} \sum_{m=0}^{N-1} a_m  \leq \liminf_{\beta \uparrow 1} (1 - \beta) \sum_{m=0}^{\infty} \beta^m a_m \nonumber \\
&\leq \limsup_{\beta \uparrow 1} (1 - \beta) \sum_{m=0}^{\infty} \beta^m a_m \leq \limsup_{N \to \infty} {1 \over N} \sum_{m=0}^{N-1} a_m
\end{align}

Thus, for every strategy pair $\sigma^1, \sigma^2$, and $\epsilon > 0$, there exists $\delta_{\epsilon}$ (depending possibly on the strategies) so that
\begin{align}
& \mathbb{E}_{\sigma^1, \sigma^2}^{\mu_0} (1 - \delta_{\epsilon}) \bigg[\sum_{m=0}^{\infty} \beta_{\epsilon}^m u^1(a^1_m,a^2_m)\bigg] + \epsilon \nonumber \\
& \geq  \liminf_{N \to \infty} {1 \over N} \mathbb{E}_{\sigma^1, \sigma^2}^{\mu_0} \bigg[ \sum_{m=0}^{N-1} u^1(a^1_m,a^2_m)\bigg]  \nonumber
\end{align}

Now, let $\sigma^1_n, \sigma^2_n$ be a sequence of strategies which converge to the supremum for the average payoff. Let $\tilde{\sigma}^1_n$, $\tilde{\sigma}^2_n$ be one which comes within $\epsilon/2$ of the supremum so that
\begin{align*}
& \sup_{\sigma^1,\sigma^2} \liminf_{N \to \infty} {1 \over N} \mathbb{E}_{\sigma^1,\sigma^2} \bigg[ \sum_{m=0}^{N-1} u^1(a^1_m,a^2_m)\bigg] \\
& \leq \liminf_{N \to \infty} {1 \over N} \mathbb{E}_{\tilde{\sigma}^1_n,\tilde{\sigma}^2_n} \bigg[ \sum_{m=0}^{N-1} u^1(a^1_m,a^2_m)\bigg] + \epsilon/2
\end{align*}

Let now $\delta_{\epsilon}$ close to 1 be a discount factor whose optimal payoff comes within $\epsilon/2$ of the limit when $\delta=1$. For this parameter, under $\tilde{\sigma}^1_n,\tilde{\sigma}^2_n$ one obtains an upper bound on this payoff, which can be further upper bounded by optimizing over all possible strategies for this $\delta_{\epsilon}$ value. This leads to a stationary strategy. Thus,
\begin{eqnarray}
&& \sup_{\sigma^1,\sigma^2} \liminf_{N \to \infty} {1 \over N} \mathbb{E}_{\sigma^1,\sigma^2} \bigg[ \sum_{m=0}^{N-1} u^1(a^1_m,a^2_m)\bigg] -\epsilon/2 \nonumber \\
&& \leq  \liminf_{N \to \infty} {1 \over N} \mathbb{E}_{\tilde{\sigma}^1_n,\tilde{\sigma}^2_n} \bigg[ \sum_{m=0}^{N-1} u^1(a^1_m,a^2_m)\bigg] \bigg] \nonumber\\
&& \leq  \mathbb{E}_{\tilde{\sigma}^1_n,\tilde{\sigma}^2_n} (1 - \delta_{\epsilon}) \bigg[\sum_{m=0}^{\infty} \delta_{\epsilon}^m u^1(a^1_m,a^2_m)\bigg] + \epsilon/2 \nonumber \\
&& \leq \mathbb{E}_{\tilde{\sigma}^1,\tilde{\sigma}^2}(1 - \delta_{\epsilon}) \bigg[\sum_{m=0}^{\infty} \delta_{\epsilon}^m u^1(a^1_m,a^2_m)\bigg] + \epsilon/2  \nonumber \\
&&\leq  \limsup_{N \to \infty} {1 \over N}  \mathbb{E}_{\tilde{\sigma}^1,\tilde{\sigma}^2} \bigg[ \sum_{m=0}^{N-1} u^1(a^1_m,a^2_m)\bigg] + \epsilon/2 + \epsilon'  \label{bound'}\\
&& =  \liminf_{N \to \infty} {1 \over N}  \mathbb{E}_{\tilde{\sigma}^1,\tilde{\sigma}^2} \bigg[ \sum_{m=0}^{N-1} u^1(a^1_m,a^2_m)\bigg] + \epsilon/2 + \epsilon'  \label{bound''}
\end{eqnarray}
where $\epsilon'$ in (\ref{bound'}) is a consequence of the following analysis. Under any stationary optimal strategy $\tilde{\sigma}^1,\tilde{\sigma}^2$ for a discounted problem,
\begin{align}
&\mathbb{E}_{\tilde{\sigma}^1,\tilde{\sigma}^2}(1 - \delta_{\epsilon}) \bigg[\sum_{m=0}^{\infty} \delta_{\epsilon}^m u^1(a^1_m,a^2_m)\bigg] \nonumber \\
&-  \limsup_{N \to \infty} {1 \over N}  \mathbb{E}_{\tilde{\sigma}^1,\tilde{\sigma}^2} \bigg[ \sum_{m=0}^{N-1} u^1(a^1_m,a^2_m)\bigg] \label{unifS}
\end{align}
is uniformly bounded over all stationary policies under Assumption \ref{unifIdentifiability}. Note finally that since $\tilde{\sigma}^1$ is stationary, limit infimum in (\ref{bound''}) and limit supremum in (\ref{bound'}) are identical by an application of the dominated convergence theorem (since the actual limit exists as $N \to \infty$). Thus, one can select $\epsilon'$ and then $\epsilon$ arbitrarily small so that the result holds in the following fashion: First pick $\epsilon' > 0$, find a corresponding $\delta_{\epsilon'}$ with the understanding that for all $\delta_{\epsilon} \in [\delta_{\epsilon'},1)$, (\ref{bound'}) holds. Now select $\delta_{\epsilon} \geq \delta_{\epsilon'}$ to satisfy the second inequality, such a $\delta_{\epsilon}$ is guaranteed to exist since there are infinitely many such $\delta$ values up to $1$ that satisfies this inequality. Here the uniformity of the convergence in (\ref{unifS}) over all stationary policies is crucial.

In the above analysis, $\tilde{\sigma}^1, \tilde{\sigma}^2$ are stationary and with this stationary strategy,
\[\lim_{N \to \infty} {1 \over N} \mathbb{E}^{\mu^1, \mu^2}_{\mu_0}[\sum_{m=0}^{N-1} u^1(a^1_m,a^2_m)] \to \int \nu^*(d\mu,\gamma) G(\mu,\gamma) \]
by the convergence of the expected empirical occupation measures, where $\nu^*$ is some invariant probability measure induced by some optimal stationary strategy. Observe also that such an optimal stationary strategy places a dirac delta measure on the normal type given the stated observability assumptions under its invariant probability measure (which in turn is a stage-wise commitment policy). This leads to the following result which says that the supremum over all strategies is equal to the supremum over stationary strategies which satisfy the structure given in Lemma \ref{walrandType}, let us call such strategies $\Sigma_M$ :
\begin{eqnarray}
&& \sup_{\sigma^1} \liminf_{N \to \infty} {1 \over N} \mathbb{E}_{\sigma^1, \sigma^2}^{\mu_0} \sum_{m=0}^{N-1} u^1(a^1_m,a^2_m)  \nonumber \\
&& = \sup_{\sigma^1 \in \Sigma_M} \liminf_{N \to \infty} {1 \over N} \mathbb{E}_{\sigma^1, \sigma^2}^{\mu_0 = \nu^*} \sum_{m=0}^{N-1} u^1(a^1_m,a^2_m)
\end{eqnarray}
Accordingly by Assumption \ref{observabilityC}, the invariant measure on $\mu_t$ will place a full mass on this type and by (\ref{whenInvariant}), we conclude that an optimal strategy exists for Player 1, which will be of commitment type. Finally, we establish that this payoff is attainable for an arbitrary initial prior satisfying the stated assumptions:
\begin{eqnarray}
&& \sup_{\sigma^1} \liminf_{N \to \infty} {1 \over N} \mathbb{E}_{\sigma^1, \sigma^2}^{\mu_0} \sum_{m=0}^{N-1} u^1(a^1_m,a^2_m)  \nonumber \\
&& = \sup_{\sigma^1 \in \Sigma_M} \liminf_{N \to \infty} {1 \over N} \mathbb{E}_{\sigma^1, \sigma^2}^{\mu_0} \sum_{m=0}^{N-1} u^1(a^1_m,a^2_m)
\end{eqnarray}
This follows from the fact that,
\begin{eqnarray}
&& \sup_{\sigma^1} \liminf_{N \to \infty} {1 \over N} \mathbb{E}_{\sigma^1, \sigma^2}^{\mu_0} \sum_{m=0}^{N-1} u^1(a^1_m,a^2_m)  \nonumber \\
&& \geq \sup_{\sigma^1 \in \Sigma_M} \liminf_{N \to \infty} {1 \over N} \mathbb{E}_{\sigma^1, \sigma^2}^{\mu_0 = \mu^*} \sum_{m=0}^{N-1} u^1(a^1_m,a^2_m)
\end{eqnarray}
and that by the identifiability condition Assumption \ref{observabilityC}, the same expected payoff (induced by the Stackelberg mimicking commitment strategy) is incurred for every initial prior (satisfying the aforementioned absolute continuity condition; that is, the full-support prior condition)%and continuity condition URBANA} through using the Stackelberg strategies optimal for $\mu_0 = \mu^*$ to an arbitrary $\mu_0$, one obtains that
\begin{eqnarray}
&& \inf_{\sigma^1} \liminf_{N \to \infty} {1 \over N} \mathbb{E}_{\sigma^1, \sigma^2}^{\mu_0} \sum_{m=0}^{N-1} u^1(a^1_m,a^2_m)  \nonumber \\
&& - \inf_{\sigma^1 \in \Sigma_M} \liminf_{N \to \infty} {1 \over N} \mathbb{E}_{\sigma^1, \sigma^2}^{\mu_0 = \mu^*} \sum_{m=0}^{N-1} u^1(a^1_m,a^2_m)
\nonumber  \\
&& = 0
\end{eqnarray}
Thus any optimal strategy will need to be infinite repetition of a stage game Stackelberg action. \qed 

\subsection{Proof of Lemma \ref{weakContKernel}.}
From (\ref{walvarRecursion}), we observe the following. Let $f$ be a continuous function on $\Delta(\Omega)$. Then $E[f(\mu_{t+1})|\mu_t,\gamma^1_t]$ is continuous in $(\mu_t, \gamma^1_t)$ if
\[\sum_{s^2_t} f(H(\mu_t,s^2_t,\gamma^1_t)) P_\sigma(s^2_t|\gamma^1_t)\] is continuous in $\mu_t, \gamma^1_t$  where $\mu_{t+1} = H(\mu_t,s^2_t,\gamma^1_t)$ defined by (\ref{walvarRecursion}) with the variables
\begin{equation}
\resizebox{0.99\hsize}{!}{$ 1_{\{\gamma^1_t(\omega,s^2_{[0,t-1]}) = a^1_{t}\}} = P_\sigma(a^1_{t} | \omega, s^2_{[0,t-1]}), \quad \mu_t(\omega)= P_\sigma(\omega |s^2_{[0,t-1]})$} \nonumber
\end{equation}
Instead of considering continuous functions on $\Delta(\Omega)$, we can also consider continuity of $\mu_{t+1}(\omega)$ for every $\omega$ since pointwise convergence implies convergence in total variation by Scheff\'{e}'s Theorem, which in turn implies weak convergence. Now, for every fixed $s^2_t=s$, $\mu_{t+1}(\omega)$ is continuous in $\mu_t$ for every $\omega$, and hence $H(\mu_t,s^2_t,\gamma^1_t)$ is continuous in total variation since pointwise convergence implies convergence in total variation. Furthermore, $P_\sigma(s^2_t|\gamma^1_t,\mu_t)$ is continuous in $\mu_t$ for a given $\gamma^1_t$; thus, weak continuity follows. 
\qed

\subsection{Proof of Lemma \ref{boundLearning}.}
%\sy{We assume that the commitment type plays $a^1$ always -needs to be written in the statement also.-}
{Suppose that $\max_x u^2(a^1,x)=u^2(a^1,x^*)$. Let \\ $P_\sigma(a^1|s^2_{[0,t]}) \geq 1 - \epsilon$. Let the maximum of
\[P_\sigma(a^1|s^2_{[0,t]}) u^2(a^1,x) + \sum_{\bar{a}^1_j \neq a^1} P_\sigma(\bar{a}^1_j|s^2_{[0,t]}) u^2(\bar{a}^1_j, x)\]
be achieved by $x^*$ so that
\begin{eqnarray}
 P_\sigma(a^1|s^2_{[0,t]}) u^2(a^1,x') +\sum_{\bar{a}^1_j \neq a^1} P_\sigma(\bar{a}^1_j|s^2_{[0,t]})u^2(\bar{a}^1_j, x')  \nonumber \\ 
\leq P_\sigma(a^1|s^2_{[0,t]}) u^2(a^1,x^*) + \sum_{\bar{a}^1_j \neq a^1} P_\sigma(\bar{a}^1_j|s^2_{[0,t]}) u^2(\bar{a}^1_j,x^*) \nonumber
\end{eqnarray}
for any $x'$. For this to hold it suffices that
\[ P_\sigma(a^1|s^2_{[0,t]}) (u^2(a^1,x^*) - u^2(a^1,x')) \geq \max_{s,t} \epsilon u^2(s,t)  \]
and since $P_\sigma(a^1|s^2_{[0,t]})  \geq 1 - \epsilon$,
\[  (u^2(a^1,x^*) - u^2(a^1,x')) \geq {\max_{s,t} \epsilon u^2(s,t)  \over 1-\epsilon}. \]
Thus, if $P_\sigma(a^1|s^2_{[0,t]}) > \epsilon$ then the optimal response is to $a^1$. In particular, with $P_\sigma(a^1|s^2_{[0,t]})  \geq 1 - \epsilon$ and for all $\bar{a}^1_j \neq a^1$ we have $P_\sigma(\bar{a}^1_j|s^2_{[0,t]}) \leq \epsilon/M$, (\ref{criterion}) holds.}
\qed

\subsection{Proof of Theorem \ref{geometricBound}.}
(\ref{criterion}) is equivalent to, by Bayes' rule:
\[{P_\sigma(s^2_{[0,t]} | \hat{\omega} = m) \over P_\sigma(s^2_{[0,t]} | \hat{\omega} = k)} \geq {P_\sigma(\hat{\omega}=k) f(M)  \over P_\sigma(\hat{\omega}=m)}\]
and
\[\sum_{j=0}^n \log{(P_\sigma(s^2_j|\hat{\omega}=m) \over (P_\sigma(s^2_j|\hat{\omega}=k)} \geq \log\bigg( {P_\sigma(\hat{\omega}=k) f(M) \over P_\sigma(\hat{\omega}=m)} \bigg)\]
Note now that (\ref{criterion}) implies that $t \subset \tau_m$. Thus, we can now apply a measure concentration result through McDiarmid's inequality (see \cite{raginsky2012concentration}) to deduce that
\begin{eqnarray}
&& P_\sigma(t \notin \tau_m) \nonumber \\
&&\resizebox{0.9\hsize}{!}{$ \leq P\bigg( \sum_{j=0}^t \log({P_\sigma(s^2_j|\hat{\omega}=m) \over P_\sigma(s^2_j|\hat{\omega}=k)}) \leq \log({P_\sigma(\hat{\omega}=k) f(M) \over P_\sigma(\hat{\omega}=m)}) \bigg)  \nonumber$} \\
&&\resizebox{0.9\hsize}{!}{$ \leq P\bigg({1 \over t+1} \sum_{j=0}^t \log({P_\sigma(s^2_j|\hat{\omega}=m) \over P_\sigma(s^2_j|\hat{\omega}=k)}) - \mathbb{E}[\log{(P_\sigma(s^2_j|\hat{\omega}=m) \over (P_\sigma(s^2_j|\hat{\omega}=k)}]    \nonumber$} \\
&& \resizebox{0.9\hsize}{!}{$ \quad \quad \quad \quad \leq {1 \over t+1} \log({P_\sigma(\hat{\omega}=k) f(M) \over P_\sigma(\hat{\omega}=m)}) - \mathbb{E}[\log{(P_\sigma(s^2_j|\hat{\omega}=m) \over (P_\sigma(s^2_j|\hat{\omega}=k)}]  \bigg)  \nonumber $}\\
&& \resizebox{0.9\hsize}{!}{$ \leq P\bigg( \bigg| {1 \over t+1} \sum_{j=0}^t \log({P_\sigma(s^2_j|\hat{\omega}=m) \over P_\sigma(s^2_j|\hat{\omega}=k)}) - \mathbb{E}[\log{(P_\sigma(s^2_j|\hat{\omega}=m) \over (P_\sigma(s^2_j|\hat{\omega}=k)}] \bigg| \nonumber $}\\
&& \resizebox{0.9\hsize}{!}{$ \quad \quad \quad \quad  \geq  | \mathbb{E}[\log{(P_\sigma(s^2_j|\hat{\omega}=m) \over (P_\sigma(s^2_j|\hat{\omega}=k)}]  - {1 \over t+1} \log({P_\sigma(\hat{\omega}=k) f(M) \over P_\sigma(\hat{\omega}=m)}) | \bigg)  \nonumber$} \\
&& \leq 2 e^{-t \bigg(\mathbb{E}[\log{(P_\sigma(s^2_j|\hat{\omega}=m) \over (P_\sigma(s^2_j|\hat{\omega}=k)}]  - {1 \over t+1} \log({P_\sigma(\hat{\omega}=k) f(M) \over P_\sigma(\hat{\omega}=m)}) \bigg)^2/(b-a)}
\end{eqnarray}
where $a \leq \mathbb{S}^j \leq b$ with $\mathbb{S}^j = {P_\sigma(s^2_j|\hat{\omega}=m) \over P_\sigma(s^2_j|\hat{\omega}=k)}$.
This implies that
%, since $\log(n)/n \to 0$ and $f(n) = n (1-\epsilon)/\epsilon$ and by Lemma \ref{boundLearning}, 
the probability of $t \notin \tau_m$ is upper bounded asymptotically by a geometric random variable, that is, there exists $R < \infty$ and $\rho \in (0,1)$ so that for all $t \in \mathbb{N}$, $P_\sigma(t \notin \tau_{m}) \leq R \rho^{t}$. \hfill \qed

%\begin{IEEEbiography} {\bf  Nuh Ayg\"un Dalk\i ran} received his B.Sc. degrees in Economics and Mathematics from Middle East Technical University in 2004, his M.A. degree in Economics from Sabanc\i \ Univesity in 2006, and his Ph.D. degree in Managerial Economics and Strategy from Kellogg School of Management, Northwestern University in 2012. He started working as an Assistant Professor in the Department of Economics at Bilkent University in 2012. His research interests are on game theory, repeated games and reputations, economics of information, and decision theory.
%\end{IEEEbiography}
%
%\begin{IEEEbiography}
%{\bf Serdar Y\"uksel} (S'02, M'11) received his B.Sc. degree in Electrical and Electronics
%Engineering from Bilkent University; M.S. and Ph.D. degrees in
%Electrical and Computer Engineering from the University of Illinois at Urbana-
%Champaign in 2003 and 2006, respectively. He was a post-doctoral researcher
%at Yale University before joining the Department of Mathematics and Statistics at Queen's University. His research interests are on stochastic control, decentralized control, information theory, and probability. He has been an Associate Editor for the IEEE TRANSACTIONS ON AUTOMATIC CONTROL, Automatica, and Systems and Control Letters.
%\end{IEEEbiography}

\end{document}